\documentclass[a4paper]{amsart} 
\usepackage[ascii]{inputenc} 

\usepackage{amssymb}
\usepackage{mathrsfs}
\usepackage[hidelinks]{hyperref}

\usepackage{xcolor}

\usepackage[shortlabels]{enumitem}
\setlist[enumerate,1]{label={(\Alph*)}}
\setlist[enumerate,2]{label={(\alph*)}}
\setlist[enumerate,3]{label={$\bullet_{\arabic*}$}}

\newenvironment{PROOF}[2][\proofname.]
{\begin{proof}[#1]\renewcommand{\qedsymbol}{\textsquare$_{\rm #2}$}}
{\end{proof}}

\newtheorem{theorem}{Theorem}[section] 
\newtheorem{claim}[theorem]{Claim}

\newtheorem{conclusion}[theorem]{Conclusion}
\newtheorem{observation}[theorem]{Observation}

\theoremstyle{definition}
\newtheorem{convention}[theorem]{Convention}
\newtheorem{definition}[theorem]{Definition}
\newtheorem{example}[theorem]{Example}
\newtheorem{problem}[theorem]{Problem}
\newtheorem{fact}[theorem]{Fact}

\newtheorem{discussion}[theorem]{Discussion}

\newtheorem{thesis}[theorem]{Thesis}

\theoremstyle{remark}

\newtheorem{remark}[theorem]{Remark}
\newtheorem{question}[theorem]{Question}
\newtheorem{notation}[theorem]{Notation}

\newtheorem{conjecture}[theorem]{Conjecture}

\newcommand{\Random}{\mathsf{Random}}
\newcommand{\Sacks}{\mathsf{Sacks}}

\newcommand{\EF}{\mathrm{EF}}
\newcommand{\GEM}{\mathrm{GEM}}

\newcommand{\QFT}{\mathrm{QFT}}

\newcommand{\cof}{\mathsf{cof}}
\newcommand{\fr}{\mathrm{fr}}

\newcommand{\oi}{\mathrm{oi}}

\newcommand{\ptr}{\mathrm{ptr}}
\newcommand{\str}{\mathrm{str}}

\newcommand{\dorg}{\mathrm{dorg}}
\newcommand{\org}{\mathrm{org}}
\newcommand{\x}{\mathrm{x}}

\newcommand{\Pos}{\mathsf{pos}}
\newcommand{\Neg}{\mathsf{neg}}

\newcommand{\iso}{\mathrm{iso}}
\newcommand{\spec}{\mathrm{spec}}

\newcommand{\Leb}{\mathrm{Leb}}
\newcommand{\lev}{\mathrm{lev}}

\DeclareMathAlphabet{\mathtxfrak}{OT1}{tx-frak}{m}{n}

\DeclareMathAlphabet{\matheulerfrak}{U}{euf}{m}{n}
\newcommand{\gk}{\matheulerfrak{k}}

\newcommand{\iif}{\mathsf{if}}
\newcommand{\oor}{\mathrm{or}}

\newcommand{\nll}{\mathrm{null}}

\newcommand{\tthh}{\mathrm{th}}

\newcommand{\NU}{\mathrm{nu}}

\newcommand{\arity}{\mathrm{arity}}

\newcommand{\cf}{\mathrm{cf}}

\newcommand{\cov}{\mathrm{cov}}

\newcommand{\dom}{\mathrm{dom}}

\newcommand{\id}{\mathrm{id}}

\newcommand{\inv}{\mathrm{inv}}

\newcommand{\Ord}{\mathrm{Ord}}

\newcommand{\rang}{\mathrm{rang}}

\newcommand{\qf}{\mathrm{qf}}

\newcommand{\CH}{\mathsf{CH}}
\newcommand{\GCH}{\mathsf{GCH}}

\newcommand{\ZFC}{\mathsf{ZFC}}

\newcommand{\Cohen}{\mathsf{Cohen}}

\newcommand{\Card}{\mathrm{Card}}

\newcommand{\EC}{\mathrm{EC}}

\newcommand{\lex}{\mathrm{lex}}

\newcommand{\LST}{\mathrm{LST}}

\newcommand{\meagre}{\mathrm{meagre}}

\newcommand{\Mod}{\mathrm{Mod}}

\newcommand{\set}{\mathrm{set}}

\newcommand{\tp}{\mathrm{tp}}

\newcommand{\tr}{\mathrm{tr}}

\newcommand{\PC}{\mathrm{PC}}

\newcommand{\Reg}{\mathrm{Reg}}

\newcommand{\bfF}{\mathbf{F}}
\newcommand{\bfG}{\mathbf{G}}

\newcommand{\bfI}{\mathbf{I}}

\newcommand{\bfL}{\mathbf{L}}

\newcommand{\bfS}{\mathbf{S}}
\newcommand{\bfT}{\mathbf{T}}

\newcommand{\bfV}{\mathbf{V}}

\newcommand{\bfX}{\mathbf{X}}

\newcommand{\bfa}{\mathbf{a}}

\newcommand{\bfc}{\mathbf{c}}

\newcommand{\bfm}{\mathbf{m}}
\newcommand{\bfn}{\mathbf{n}}
\newcommand{\bfo}{\mathbf{o}}
\newcommand{\bfp}{\mathbf{p}}

\newcommand{\bbE}{\mathbb{E}}   

\newcommand{\bbL}{\mathbb{L}}

\newcommand{\bbP}{\mathbb{P}}
\newcommand{\bbQ}{\mathbb{Q}}

\newcommand{\mn}{\medskip\noindent}
\newcommand{\sn}{\smallskip\noindent}
\newcommand{\bn}{\bigskip\noindent}

\newcommand{\cL}{\mathscr{L}}

\newcommand{\clB}{\mathcal{B}}

\newcommand{\clH}{\mathcal{H}}

\newcommand{\clM}{\mathcal{M}}

\newcommand{\clS}{\mathcal{S}}
\newcommand{\clT}{\mathcal{T}}
\newcommand{\clU}{\mathcal{U}}

\newcommand{\clW}{\mathcal{W}}
\newcommand{\clX}{\mathcal{X}}

\newcommand{\gC}{\mathfrak{C}}

\newcommand{\eps}{\varepsilon}
\newcommand{\cl}{c\kern-.11ex \ell}
\newcommand{\lh}{{\ell\kern-.27ex g}}
\newcommand{\rest}{\restriction}

\newcommand{\caret}{{\char 94}}
\newcommand{\LL}{\langle}
\newcommand{\RR}{\rangle}

\newcommand{\subref}[1]{$_{\mathrm{\texttt{=}}\mathsf{L{#1}}}$}

\newcommand{\lepref}[1]{({<}\,{#1})}

\newcommand{\overbar}[1]{\mkern 1.5mu\overline{\mkern-1.5mu#1\mkern-1.5mu}\mkern 1.5mu}

\newcommand{\ol}{\overline}
\newcommand{\olsi}[1]{\,\overline{\!{#1}}} 

\newcommand*{\defeq}{\mathrel{\vcenter{\baselineskip0.5ex \lineskiplimit0pt\hbox{\scriptsize.}\hbox{\scriptsize.}}}=}

\usepackage{tcolorbox}

\newcount\skewfactor
\def\mathunderaccent#1#2 {\let\theaccent#1\skewfactor#2
\mathpalette\putaccentunder}
\def\putaccentunder#1#2{\oalign{$#1#2$\crcr\hidewidth
\vbox to.2ex{\hbox{$#1\skew\skewfactor\theaccent{}$}\vss}\hidewidth}}
\def\name{\mathunderaccent\tilde-3 }
\def\Name{\mathunderaccent\widetilde-3 }

\newbox\noforkbox \newdimen\forklinewidth
\setbox0\hbox{$\textstyle\bigcup$}
\setbox1\hbox to \wd0{\hfil\vrule width \forklinewidth depth \dp0
                        height \ht0 \hfil}
\setbox\noforkbox\hbox{\box1\box0\relax}
\def\unionstick{\mathop{\copy\noforkbox}\limits}
\def\nonfork#1#2_#3{#1\unionstick_{\textstyle #3}#2}
\def\nonforkin#1#2_#3^#4{#1\unionstick_{\textstyle #3}^{\textstyle
    #4}#2}
\setbox0\hbox{$\textstyle\bigcup$}
\setbox1\hbox to \wd0{\hfil{\sl /\/}\hfil}
\setbox2\hbox to \wd0{\hfil\vrule height \ht0 depth \dp0 width
                                \forklinewidth\hfil}
\newbox\doesforkbox
\setbox\doesforkbox\hbox{\box1\box0\relax}
\def\nunionstick{\mathop{\copy\doesforkbox}\limits}

\def\fork#1#2_#3{#1\nunionstick_{\textstyle #3}#2}
\def\forkin#1#2_#3^#4{#1\nunionstick_{\textstyle #3}^{\textstyle
    #4}#2}

\newcommand{\stickT}{%
\setbox255=\hbox{\raise1ex\hbox{$\hspace{0.2pt}\,\bullet\,$}}
\mathord{\rlap{\hbox to\wd255{\hss\hbox{$|$}\hss}}
\box255}
}
\newcommand{\stickS}{%
\setbox255=\hbox{\raise0.6ex\hbox{$\scriptstyle\bullet$}}
\mathord{\rlap{\hbox to\wd255{\hss\hbox{$\scriptstyle|$}\hss}}
\box255}
}

\author[S. Shelah]{Saharon Shelah}
\address{Einstein Institute of Mathematics,
The Hebrew University of Jerusalem,
9190401, Jerusalem, Israel; and\\
Department of Mathematics,
Rutgers University,
Piscataway, NJ 08854-8019, USA}
\urladdr{https://shelah.logic.at/}
\thanks{First typed 2024-11-12.
The author would like to thank Craig Falls for generously funding typing services, and thanks Matt Grimes for the careful and beautiful typing.\\
The author would also like to thank the Israel Science Foundation for
partial support of this research by grant 2320/23 (2023-2027).\\
References like [Sh:950, Th0.2\subref{y5}] mean that the internal label of Theorem 0.2 in Sh:950 is \textsf{y5}.
The reader should note that the version in my website is usually more up-to-date than the one in arXiv.
This is publication number    
1261
in Saharon Shelah's list.
}

\makeatletter
\@namedef{subjclassname@2020}{\textup{2020} Mathematics Subject Classification}
\makeatother
\date{June 16, 2025} 

\title[$\bbP$-twins I]{Twins: non-isomorphic models forced to be isomorphic\\ Part I --- 1261}

\begin{document}
\makeatletter\def\shfiuwefootnote{\gdef\@thefnmark{}\@footnotetext}\makeatother\shfiuwefootnote{Version 2025-06-17\_2. See \url{https://shelah.logic.at/papers/1261/} for possible updates.}
\begin{abstract}
    For which (first-order complete, usually countable) $T$ do there exist non-isomorphic models of $T$ which become isomorphic after forcing with a forcing notion $\bbP$? Necessarily, $\bbP$ is non-trivial; i.e.~it adds some new set of ordinals. It is best if we also demand that it collapses no cardinal. It is better if we demand on the one hand that the models are non-isomorphic, and even \emph{far} from each other (in a suitable sense), but on the other hand, $\cL$-equivalent in some suitable logic $\cL$. 
    
    In this part we give sufficient conditions: for theories with the independence property, we prove this when $\bbP$ adds no new $\omega$-sequence. We may prove it ``for some $\bbP$," but better would be for some specific forcing notions, or a natural family. Best would be to characterize the pairs $(T,\bbP)$ {for which we have such models}.

    The results say (e.g.) that there are models $M_1,M_2$ which are not isomorphic (and even \emph{far} from being isomorphic, in a rigorous sense) which become isomorphic when we extend the universe by adding a new branch to the tree $({}^{\theta>}2,\lhd)$. 

    We shall mention some specific choices of $\bbP$: mainly 
    $({}^{\theta>}2,\lhd)$ with $\theta = \theta^{<\theta}$.

    This work does not require any serious knowledge of forcings, nor of stability theory, 
though they form the motivation.
Concerning forcing,  
    the reader just has to agree that starting with a universe $\bfV$ of set theory (i.e.~a model of $\ZFC$) and a quasiorder 
    $\bbP$, there are a new directed $\bfG \subseteq \bbP$ meeting every dense subset $D$ of $\bbP$ and a universe $\bfV[\bfG]$ (so {it} satisfies $\ZFC$) of which the original $\bfV$ is a transitive subclass. We may say that $\bfV[\bfG]$ (also denoted $\bfV^\bbP$) is the universe obtained by forcing with $\bbP$. 
    
    This is part of the classification and so-called \emph{Main Gap} programs.
\end{abstract}

\maketitle

\newpage

\centerline{\textbf{\underline{Annotated Content}}}

\bn 
\textbf{\S0\quad Introduction} \hfill pg.\pageref{S0}

\bn 
\textbf{\S1\quad GEM Models} \hfill pg.\pageref{S1}
\begin{enumerate}
    \item [\ref{a5} --] {Defining} $\GEM(I,\Phi)$ for a general class $K$
\sn
    \item [\ref{a8} --] The basic examples of $K$-s
\sn
    \item[\ref{a41} --] $\Upsilon_K[T,\kappa]$; blueprints of GEM models
\sn
    \item[\ref{a44} --] $\varphi$ witnesses properties of $T$
\sn
    \item[\ref{a51} --] Proving the existence of the relevant $\Phi \in \Upsilon_K[T,\kappa]$ for $\varphi$
\sn
    \item[\ref{a56} --] $\Phi$ represents $(\varphi,R)$
\sn
    \item[\ref{a59}-\ref{a68} --] Definitions of being `far'
\end{enumerate}

\bn 
\textbf{\S2\quad Toward non-isomorphic twins} \hfill pg.\pageref{S2}

\noindent
\textbf{\S2(A)}\quad The Frame
\begin{enumerate}
    \item [\ref{b3} --] Defining twinship parameters
\sn
    \item [\ref{b20} --] Defining $\Omega_\bfp$ and $F_{\eta,\iota}$
\sn
    \item [\ref{b19} --] Main definition: $K_{\clT,\ell}^\oor$ and $K_{\clT,\ell}^\org$ (for $\ell = 1,2$)
\end{enumerate}

\sn
\textbf{\S2(B)} \quad Examples

Here we derive some twinship parameters from forcing notions, and show that for some forcing notions, there is no twinship parameter $\bfp$ with $\theta_\bfp = \aleph_0$.

\sn
\textbf{\S2(C)} \quad Are the $K$-s reasonable?

\begin{enumerate}
    \item[\ref{b26} --] $K$ with JEP and amalgamation
\sn
    \item[\ref{b29} --] Proving that $K_{\clT,\ell}^\oor$ and $K_{\clT,\ell}^\org$ are AECs and are universal
\end{enumerate}

\bn 
\textbf{\S3\quad Existence for independent $T$} \hfill pg.\pageref{S3}

    Here we deal with theories with the independence property and $\theta_\bfp > \aleph_0$, and get a version of our theorems which are quite strong. 
    That is, we have a $\lepref{\theta}$-complete forcing notion $\bbP$ (or just one which adds no new sequences of ordinals of length $<\theta$). We build $\bbP$-twins of cardinality $\lambda$ (with $\lambda > 2^{\|\bbP\|}$ a regular cardinal). Moreover, these twins are far from each other (and more).

\begin{enumerate}
    \item [\ref{c5} --] Defining entangledness 
\sn
    \item [\ref{c8} --] unembeddability
\sn
    \item[\ref{c11} --] Sufficient conditions for $I_s \in K_\clT^\org$; $(\mu,\kappa)$-unembeddable into $I_c$
\sn
    \item[\ref{c14} --] Improving on Claim \ref{c11}
\sn
    \item[\ref{c17} --] Every independent $T$ is witnessed by a $\theta$-complete $\theta^+$-cc forcing notion, where $\theta = \theta^{< \theta}$.
\end{enumerate}

\newpage
\setcounter{section}{-1}

\section{Introduction}\label{S0}

We are interested in classifying theories (or classes of models --- i.e.~structures) by the possible existence of models which are very similar but not isomorphic.

\mn
\begin{definition}\label{x1}
1) For a forcing notion $\bbP$, we say the models $M$ and $N$ are $\bbP$-\emph{isomorphic} \underline{when} they become isomorphic after forcing with $\bbP$.

\mn
2) For $\bfX$ a set or class of forcing notions, we say $M$ and $N$ are \emph{strongly} $\bfX$-{isomorphic} \underline{when} they are $\bbP$-isomorphic for {every} $\bbP \in \bfX$. 

\mn
3) \emph{Weakly} $\bfX$-{isomorphic} (or simply `$\bfX$-{isomorphic}') will mean``for some $\bbP \in \bfX$." E.g.~`weakly ccc-isomorphic' means ``for {some} ccc forcing notion."
\end{definition}

\mn
\begin{definition}\label{x2.0}
1) We say two models are $\bbP$-\emph{twins} when they are $\bbP$-isomorphic but not isomorphic. We say they are $(\bbP,\cL)$-twins when they are $\bbP$-twins and $\cL$-equivalent, for $\cL$ a logic.

\mn
2)  We say $M$ and $N$ are $(\bbP,\lambda)$-twins (or $(\bbP,\cL,\lambda)$-twins) when in addition, $\|M\| = \|N\| = \lambda$.

\mn
3) Similarly for $\bfX$-twins and strong $\bfX$-twins.

\mn
4) We may say a theory $T$ [or a class $K$ of models] `has $\bbP$-twins.'
\end{definition}

Baldwin-Laskowski-Shelah \cite{Sh:464} and Laskowski-Shelah \cite{Sh:518} investigated the case of weak ccc-twins (i.e.~$\bfX$ is the class of ccc forcings).\footnote{
    There, twins were called `potentially isomorphic.'
}
Lately, Farah raised a similar question, for $\bbP$ the Random Real forcing and $T$ an unstable theory.

\mn
\centerline{*\qquad*\qquad*}

We thank the referee and Shimon Garti for their helpful comments.

\bn
\subsection{A panoramic picture -- the long-range view}

\mn
\begin{thesis}[\textbf{The Classification Thesis}]\label{x3}
We would like to classify the theories $T$: naturally, at first all complete first-order (maybe countable) ones, but later try for more --- e.g.~for every AEC.

Like Janus, the thesis has two faces:
\begin{enumerate}
    \item Set theoretic test questions which will shed light on the complexity of $T$, leading to constructing `complicated' models of a theory $T$, when $T$ itself is complicated.
\sn
    \item Finding dividing lines among the family of theories, such that
    \begin{enumerate}[$\bullet_1$]
        \item Above the line, we have results as in (A).
\sn
        \item Below the line we develop structure theory, and can analyze models of $T$ to some extent. 
    \end{enumerate}
\sn
    \item The thesis is that those two sides of the program are strongly connected, because if we succeed in proving a case of the so-called \emph{main gap}, we get complementary results. So we know that the assumptions in each are the best possible, and with this aim in mind one is driven to discover inherent properties of $T$.
\sn
    \item Even if you are only interested in clause (B)\,$\bullet_2$, this thesis tells you that having (A) \emph{and} (B) in mind is a good way to advance each of them.
\sn
    \item For a given test question, a `Main Gap' theorem will describe how the theories are divided into ones with complicated models (the \emph{non-structure} side), and ones with a `structure theory.'

    But naturally, along the way we may come across other properties which could be of interest, possibly more than the original question (e.g.~whether $T$ is stable).
\sn
    \item Having the two sides gives us more than the sum of their parts; it proves that both are maximal (in the chosen context), and that those properties are the natural dividing line. 

    (Of course, not all interesting properties are like this: you may be able to say something about binary functions on a set which is not a group, but this is not the animating question on the class of groups. Closer are 
    {o}-minimal theories.)
\end{enumerate}
\end{thesis}

\bn
\centerline{*\qquad*\qquad*}

\bigskip

The classical case was first-order complete countable theories, but there are others; {e.g.} universal classes {up to AEC}.

In \cite{Sh:a} and \cite{Sh:c}, the set theoretic test questions were:
\begin{itemize}
    \item $I(\lambda,T) \defeq$ the number of isomorphism classes of models in $\EC_\lambda(T)$ (= models of $T$ of cardinality $\lambda$).
\sn
    \item $IE(\lambda,T) \defeq$ the maximal number of pairwise-non-elementarily embeddable models in $\EC_\lambda(T)$.
\end{itemize}
In this case, the thesis was that this classification characterizes answers to the question ``Is $\Mod_T$ (the class of models of $T$) complicated?", along a significant number of measures.

With regards to those test questions, the situation can be seen in the following trichotomy; the uninitiated reader may concentrate on $\boxplus_2,\boxplus_3$. 

This theorem is the original case: the Main Gap Theorem of \cite{Sh:c}.
\begin{theorem}[\textbf{The main gap Trichotomy}]\label{x17}
$1)$
\begin{enumerate}[$(A)$]
    \item [$\boxplus_1$] The countable complete first-order theories $T$ can be divided into three classes:
\sn
    \item Unstable \underline{or} stable but un\emph{super}stable \underline{or} superstable with OTOP \underline{or} superstable with DOP.
    (The last two cases tell us some non-first-order formula defining many graphs in some models of $T$.)
\sn
    \item $T$ is not in $(A)$, but it is \emph{deep}.
\sn
    \item Neither $(A)$ nor $(B)$. (The antonym of deep is \emph{shallow}.)
\end{enumerate}

\mn
$2)$ The classification in part $(1)$ is by the {inside properties of these theories; this is not meaningful if you do not know them.}

 Let us move to the other side of the coin.
\begin{enumerate}
    \item [$\boxplus_2$] If $T$ is of type $(B)$ or $(C)$ it is called \emph{classifiable}, and satisfies the following:
    \begin{enumerate}[$(a)$]
        \item A model $M$ of $T$ can be described by a tree $\clT$ with 
        $\omega$ levels. That is, it is a set of finite sequences, closed under initial segments (and countably many unary predicates).
\sn
        \item More fully, there is a tree $\LL M_\eta : \eta \in \clT\RR$ of countable submodels,\\ $\prec$-increasing with $\eta$,
        \emph{``freely joined"} (i.e.~this tree of models is non-forking), and $M$ is prime over 
        $\bigcup\limits_{\eta \in \clT} M_\eta$.
\sn
        \item Another aspect is: models of $T$ can be characterized (up to isomorphism) by their theory in the logic $\cL = \bbL_{\infty,\aleph_1}$, enriched by ``cardinality quantifiers on dimension by definable dependence relations" (see \cite[Ch.XIII]{Sh:c}, \cite{Sh:307}, \cite{Sh:897}).
    \end{enumerate}
\end{enumerate}

\mn
$3)$ Continuing in this fashion:
\begin{enumerate}
    \item [$\boxplus_3$] If $T$ satisfies $(A)$, then a strong negation of the above holds. The class of models is \emph{non-classifiable}: e.g.~models (pedantically, isomorphism classes of models) code stationary sets. Specifically, for a model $M$ of $T$ of cardinality $\lambda$ ($\lambda$ regular uncountable) we can find an invariant $\inv(M) = \inv(M/{\cong})$ of the form $S/\mathsf{club}$, for some stationary $S \subseteq \lambda \cap \cof(\aleph_0)$, {so that every such $S/\mathsf{club}$ occurs} (see \cite[2.4, 2.5(2), pp.296-7]{Sh:220}).
\end{enumerate}

\mn
$4)$ More on {$\boxplus_1$} ---
\begin{enumerate}
    \item [$\boxplus_4$] 
    \begin{enumerate}[$(a)$]
        \item If $T$ satisfies $\boxplus_1(A)$ or $(B)$, \underline{then} for every cardinal $\lambda$, $T$ has $2^\lambda$-many pairwise non-isomorphic models of cardinality $\lambda$ (the maximal number possible). 
\sn
        \item If $T$ satisfies $\boxplus_1(C)$ \underline{then}
        $I(\aleph_\alpha,T) < \beth_{\omega_1}(|\alpha|)$ 
        for every ordinal $\alpha$. So it fails the conclusion of clause $(a)$ when (e.g.) $\GCH$ holds.
\sn
        \item Suppose $T$ satisfies $\boxplus_1(C)$. If $\LL M_\alpha : \alpha < \beth_{\omega_1}\RR$ is a sequence of models of $T$, \underline{then} for some $\alpha < \beta < \beth_{\omega_1}$ there is an elementary embedding of $M_\alpha$ into $M_\beta$.
\sn
        \item Suppose $T$ satisfies $\boxplus_1(A)$. \underline{Then} for every $\lambda > \aleph_0$ there exists a family of $2^\lambda$-many models of $T$, each of cardinality $\lambda$, with no one elementarily embeddable into another.
\sn
        \item For $T$ that satisfy $\boxplus_1(B)$, their behavior is in the middle --- for some cardinal $\kappa$ (the first so-called \emph{beautiful} cardinal\footnote{
            A \emph{beautiful} cardinal is a large cardinal which is compatible with `$\bfV = \bfL$,' but whose existence cannot be proved in $\ZFC$.
        }) 
        we have:
        \begin{enumerate}
            \item If $\lambda \in (\aleph_0,\kappa)$ it behaves as in clause $(d)$ above.
\sn
            \item If $\LL M_\alpha : \alpha < \kappa\RR$ is a sequence of models of $T$ (of any cardinality), \underline{then} for some $\alpha < \beta < \kappa$ there is an elementary embedding of $M_\alpha$ into $M_\beta$.
        \end{enumerate}
        (We may say $\kappa = \infty$ when no such cardinal exists.)
    \end{enumerate}
\end{enumerate}
\end{theorem}

\mn
But of course, there are other worthwhile measures:
\begin{problem}\label{x20}
What if we ask for which $T$-s do we have a weaker version of \ref{x17}$\boxplus_2$, where we replace $\clT$ with a tree with \underline{$\omega_1$ levels}? I.e.~$\clT$ consists of  sequences of countable length: say, subtrees of $({}^{\omega_1>}\!\lambda,\lhd)$.
\end{problem}

This calls for a finer discussion of stable theories.

\bn
\centerline{*\qquad*\qquad*}

\bigskip

\subsection{First approach}

Very similar, but not the same. 

True dividing lines (and measures of complexity) discussed above are relevant for a significant set of questions which are not \emph{a priori} connected. A major case is the \emph{Keisler order}, resolved for stable $T$. (See a recent survey by Keisler \cite{Ke17} on this.)
Another measure is the number of $\aleph_1$-\emph{resplendent} models in $\EC_\lambda(T)$, up to isomorphism (see \cite{Sh:363}, which characterizes stable $T$). Still another direction is building \emph{somewhat rigid} models (see \cite{Sh:800} and references there).

There are also works on unstable theories -- simple, dependent, and NTP$_2$ -- but here we concentrate on dividing lines among stable $T$.

We may ask for the number of $|T|^+$-saturated models of $T$ (or complete metric spaces) -- see \cite{Sh:1238}. But closer to our problem is the following way to strengthen our non-structure side:

\begin{question}\label{x23}
When do there exist models of a theory $T$ (that is, an elementary class) which are very similar but \emph{not} isomorphic? This question can serve as a yardstick for the complexity of $T$, and thus makes for a good test problem.
\end{question}

One interpretation of ``$M$ and $N$ are very similar" is 
\begin{enumerate}
    \item [$\bullet$] $M$ and $N$ are of cardinality $\lambda$, and are equivalent for a `strong logic' $\cL$.
\end{enumerate}
We call this \emph{the first approach}. 

\begin{discussion}\label{x24}
We provide references for some relevant works (including those earlier ones asking a more basic question: are there such models not restricting $T$?).
\begin{enumerate}
    \item [(A)$_1$] Existence of $\bbL_{\infty,\lambda}$-equivalent but non-isomorphic models of cardinality $\lambda$:
    
    \begin{itemize}
        \item For $\lambda$ regular uncountable, this is an unpublished result of Morley.
\sn
        \item \cite{Sh:188} covers singular $\lambda = \lambda^{\aleph_0}$.
\sn
        \item \cite[Ch.II, 7.4-5, p.111]{Sh:g} proves it for almost all 
        $$
        {\lambda^{\aleph_0} > \lambda} > \cf(\lambda) > \aleph_0.
        $$
    \end{itemize}
\sn
    \item [(A)$_2$] For $M_*$ a model of cardinality $\lambda$, what can be said about the value of 
    $$
    \NU(M_*) \defeq |K_{M_*} / {\cong}|,
    $$ 
    where $K_{M_*} \defeq \{M : \|M\| = \lambda,\ M \equiv_{\bbL_{\infty,\lambda}} M_*\}$?
    \begin{enumerate}
        \item {Palyutin} \cite{Pal77}: If $\bfV = \bfL$ and $\lambda = \aleph_1$, then $\NU(M_*) \in \{1,2^{\aleph_1}\}$.
\sn
        \item By \cite{Sh:129}: if $\bfV = \bfL$ and $\lambda$ is regular uncountable but not weakly compact,  then $\NU(M_*) \in \{1,2^\lambda\}$.
\sn
        \item By \cite{Sh:125}, the `$\bfV = \bfL$' in clause (b) is necessary. (That is, it cannot be proved in $\ZFC$.)
\sn
        \item By \cite{Sh:133}, if $\lambda$ is weakly compact and $\theta \in [1,\lambda]$, \underline{then} there exists a model $M$ with cardinality $\lambda$ and $\NU(M) = \theta$. 
    \end{enumerate}
\sn
    \item [(A)$_3$]  $\bbL_{\infty,\lambda}$-equivalent but not isomorphic models, for $T$ unsuperstable; see\\ \cite{Sh:220}.
\sn
         \item [(A)$_4$]   
         Let $\cL \defeq \bbL_{\infty,\lambda}^{(\dim)}$,  where the `$(\dim)$' means that we add quantifiers saying ``$\lambda$ is the dimension of
         a definable dependence relation satisfying the Steinitz axioms (e.g. like linear dependence in vector spaces)." 
     By \cite[Ch.XIII, Th.1.4]{Sh:c}, for a (countable complete first-order) $T$ we have the following:
    \begin{enumerate}[$\bullet_1$]
        \item If $T$ satisfies $\boxplus_1$(B) or (C) of \ref{x17}(1), then any $\cL$-equivalent models of cardinality $\lambda$ are isomorphic. (See also \ref{x17}$\boxplus_2$(c).)
\sn
        \item If $T$ satisfies \ref{x17}(1)$\boxplus_1$(A) then the conclusion of $\bullet_1$ fails {badly} (see \cite{Sh:220}).
    \end{enumerate}
\end{enumerate}
To give more details, what we really have is a separation into three classes (recall \ref{x17}(4)$\boxplus_4$(e)). 
\begin{enumerate}
    \item [(B)$_1$] Game quantifier-equivalent but not isomorphic models of cardinality $\lambda$:\footnote{
        So this is stronger {than} $\bbL_{\infty,\lambda}$-equivalence.
    } 
    see \cite{Va95}. See earlier \cite{Sh:109} with Hodges;
    also \cite{Sh:836}, \cite{Sh:866} with Havlin, and \cite{Sh:907}.
\sn
    \item [(B)$_2$] For $\tau$-models $M$ and $N$, let $\EF_{\!\alpha,\lambda}(M,N)$ denote the Ehrenfeucht-Fra\"iss\'e game with $\alpha$-many moves ($\alpha$ an ordinal), each move adding $<\lambda$ elements.
    
    Like (B)$_1$ for `dividing line' $T$-s: see Hyttinen and Tuuri \cite{HyTu91}, and Hyttinen and the author \cite{Sh:474}, \cite{Sh:529}, \cite{Sh:602}.

    By \cite{HyTu91}, if $T$ is unstable and $\lambda = \lambda^{<\lambda}$, \underline{then} there are non-isomorphic models $M,N \in \EC_\lambda(T)$ which are $\Game_\zeta^\iso(M,N)$-equivalent for all $\zeta < \lambda$. (I.e.~the \textsf{ISO} player has a winning strategy: see Definition \ref{z32}(2).) 
    {This also applies to} the version using a tree $\clT \subseteq {}^{\lambda>}\!\lambda$ with no $\lambda$-branches.

    Moreover, ``$T$ has OTOP or is superstable with DOP" will suffice.
    For unsuperstable $T$ the results are weaker.

    By \cite{Sh:474}, if $T$ is a complete first-order theory which is stable but not superstable and $\lambda \defeq \mu^+$, where $\mu = \cf(\mu) \geq |T|$, \underline{then} there are $\EF_{\mu\cdot\omega,\lambda}$-equivalent but non-isomorphic models of $T$ (and even in $\PC(T_1,T)$) of cardinality $\lambda$.

    See more in \cite{Sh:529}, \cite{Sh:602}.
\sn
    \item [(C)$_1$] The present work continues \cite{Sh:897}, in a sense.
    In the second part, we intend to deal with a logic suggested there, suitable to be an analogue of \ref{x17}(3), towards \ref{x20}. We also suggest that a family of stable theories strictly containing the superstables is relevant. 
\end{enumerate} 
\end{discussion}

\bigskip

\subsection{Second approach}

The immediate impetus for this work is

\begin{conjecture}[\textbf{The Farah Conjecture}]\label{x31}
For a (first-order countable) unstable $T$, there are non-isomorphic models $M,N$ which become isomorphic when we extend the universe by adding a random real; that is, they are $\Random$-twins.
\end{conjecture}
Farah has proved this for linear orders.

\medskip
The background behind this question can be found in Baldwin-Laskowski-Shelah \cite{Sh:464} and a work with Laskowski \cite{Sh:518}. There, `similar' was defined as ``ccc-isomorphic" (see Definition \ref{x1}).

\mn
Our aim here is to try to sort this out.

\smallskip
For both approaches, a natural dream is to characterize the theories (for now, first-order complete countable) for which this occurs. The Main Gap Theorem of \cite{Sh:c} had done this for a different test question.

So Baldwin-Laskowski-Shelah \cite{Sh:464} and Laskowski-Shelah \cite{Sh:518} pose (and partly answer) the following problems.

\mn
\begin{problem}\label{x7}
\begin{itemize}
    \item Characterize the (countable) $T$ with no ccc-twins.
\sn
    \item Characterize the (countable) $T$ such that for some $T_1 \supseteq T$, `ccc-isomorphic implies isomorphic' holds in $\PC(T_1,T)$. (See \ref{z6}(5).)
\end{itemize}
\end{problem}

To explain {the choices in} \cite{Sh:464}, recall {the classification made in} \ref{x17}.
The thesis was that this classification characterizes answers to the question ``Is $\Mod_T$ (the class of models of $T$) complicated?", along a significant number of measures: e.g.~for $\dot I(\lambda,T)$ (the number of isomorphism classes of models of $T$ of cardinality $\lambda$) or $IE(\lambda,T)$ (the number of pairwise non-elementarily embeddable models of $T$ of cardinality $\lambda$).

Clearly there is a connection between the two approaches.
\begin{enumerate}
    \item [$\boxtimes_1$] If $M$ and $N$ are non-isomorphic, have different $\cL$-theories (for some logic $\cL$), and forcing with $\bbP$ preserves the $\cL$-theory of a model, \underline{then} $M$ and $N$ cannot be $\bbP$-twins.\\ (I.e.~$\not\Vdash_\bbP``M \cong N"$.)
\end{enumerate}
Now the class of ccc forcings is a natural choice, as it preserves much of what we care about; e.g.~if $\bbP$ collapses cardinals then every $T$ has (trivial) twins --- this motivated \cite{Sh:464}. 

{There it is} proved that:
\begin{enumerate}
    \item [$\boxtimes_2$]
    \begin{enumerate}
        \item If $T$ is from subclass $\boxplus_1$(A) of \ref{x17}(1), then it has ccc-twins. 
\sn
        \item However, some theories from $\boxplus_1$(C) have ccc-twins as well.
    \end{enumerate}

    \sn
     More fully (quoting \cite{Sh:518}):
\begin{quotation}
    If $T$ [is superstable and] has only countably many complete types \emph{yet} has a type of infinite multiplicity, \underline{then} there is a ccc forcing $\bbQ$ such that in any $\bbQ$-generic extension of the universe, there are non-isomorphic models $M$ and $N$ which can be forced to be isomorphic by a ccc forcing. We give examples showing that the hypothesis on the number of complete types is necessary.
\end{quotation}
\end{enumerate}

\mn
So we still do not have the answer to our questions:
\begin{question}\label{x8}
1) Can we characterize the class of $T$ which have ccc-twins?

\mn
2) Can we characterize the class of $T$ which have $\bbP$-twins,
where $\bbP$ is $\aleph_1$-complete and collapses no cardinals?
\end{question}

\mn
The motivation for \ref{x8}(2) was the following.
\begin{observation}\label{x11}
If $\bbP$ is a forcing notion which collapses no cardinal and adds no $\omega$-sequence,
\underline{then} forcing by 
$\bbP$ preserves the $\cL$-theory of a model (where $\cL$ is as in 
\emph{\ref{x24}}$(A)_4$ from \emph{\S0(B)}).
\end{observation}

\medskip
Here we mainly tried to deal with \ref{x8}(2), but after \cite{Sh:464} and \cite{Sh:518} further work was delayed. However, Farah's Conjecture gives us new inspiration to look at these questions again, for one specific ccc forcing. 

This conjecture remains open. In general, we may ask these questions for any fixed forcing notion $\bbP$: compare to \cite{Sh:464} and \cite{Sh:518}, where we asked about ``there is a ccc $\bbP$." (Instead of `$(\exists\bbP)$' {or} $\bbP$ being specified in the question, we may even try $(\forall\bbP)$.)

Note that
\begin{enumerate}
    \item [$\boxtimes_4$] 
    \begin{enumerate}
        \item If $\theta = \theta^{< \theta} > \aleph_0$ \underline{then} 
        $\Cohen_\theta \defeq ({}^{\theta>}2,\lhd)$ is a forcing as in \ref{x8}(2).
\sn
        \item Any Suslin tree $\clT$ is a ccc forcing as in \ref{x8}(2).
    \end{enumerate}
\end{enumerate}
Anyhow, if $\bbP$ is an NNS\footnote{
    \emph{NNS} means ``adds no new $\omega$-sequence."
} 
forcing then \ref{x2} answers \ref{x8}(2).

\begin{enumerate}
    \item [$\boxplus_5$] If $\theta = \theta^{<\theta} > \aleph_0$ \underline{then} after forcing with $\Cohen_\theta$, `$\theta = \theta^{<\theta} > \aleph_0$' still holds, so we have existence theorems (see below). 
\end{enumerate}
This is the approach taken in \cite{Sh:518} for the class of ccc forcing notions.

\vspace{8mm}
\subsection{The results} 

In this part, we concentrate on independent (first-order) $T$ and the forcing $\Cohen_\theta$ (for $\theta = \theta^{<\theta} > \aleph_0$), because that is where the statements and their proofs are most transparent. 

In \S1-3, we will prove

\begin{theorem}\label{x2}
$1)$ Assume $T$ is a complete countable first-order theory which has the
independence property, and $\bbP$ is $\aleph_1$-complete (or just is an NNS forcing.)
\underline{Then} $T$ has models $M$ and $N$ which are $\bbP$-twins.
(I.e.~$M \not \cong N \wedge\Vdash_\bbP ``M \cong N"$.)

\mn
$2)$ Moreover, $M$ and $N$ are \emph{far} from each other, as defined in \emph{\ref{a68}}.

\mn
$3)$ If $T_1 \supseteq T$ is also first-order countable, \underline{then} we can add ``\,$\PC(T_1,T)$ has $\bbP$-twins." (Specifically, $M$ and $N$ can be expanded to models of $T_1$.)
\end{theorem}

In Part II of this work \cite{Sh:F2433}, we shall deal with unstable $T$ (e.g.~for $\bbP \defeq \Cohen_\theta$ with $\theta = \theta^{<\theta} > \aleph_0$) and with ccc forcing notions (e.g.~Random Real) for some {of those} $T$-s. However, some stable but unsuperstable theories fail the conclusion of \ref{x2}. 

We will continue both approaches, using the relations from \cite{Sh:210} (e.g. weaker versions of entangledness). Definitions here will be phrased so as to apply to Part II as well.

\bigskip

\centerline{*\qquad*\qquad*}

\bigskip
We thank Jakob Kellner and the referee for their helpful comments, and  most of all thank M. Cardona.

\bn
\subsection{Preliminaries}

\begin{notation}\label{z5}
We will try to use standard notation.

\mn
1) $\theta,\kappa,\lambda,\mu,\chi$ will denote cardinals (infinite, if not stated otherwise). $\lambda^+$ will denote the successor of $\lambda$.

\mn
2) $\alpha,\beta,\gamma,\delta,\eps,\zeta,\xi$, $i$, and $j$ will denote ordinals. $\delta$ will be a limit ordinal unless explicitly said otherwise.

\mn
3) $k,\ell,m,n$ will denote natural numbers.

(We may abuse this somewhat and use them as indices for ordinals $< \kappa$, in statements where the default case or usual formulation is $\kappa \defeq \aleph_0$; if so, we will mention it explicitly.)

\mn
4) $\varphi, \psi$, and $\vartheta$ will be formulas; first-order, if not said otherwise.

\mn
5) For cardinals $\kappa < \lambda = \cf(\lambda)$, let 
$$
S_\kappa^\lambda \defeq \big\{\delta < \lambda : \cf(\delta) = \cf(\kappa) \big\}
$$
and
$$
S_{\leq\kappa}^\lambda \defeq \big\{\delta < \lambda : \cf(\delta) \leq \cf(\kappa) \big\}.
$$

\mn
6) $\cof(\mu)$ will denote the class of ordinals with cofinality equal to $\cf(\mu)$.

\mn
7) $\lambda^+$ may be written $\lambda(+)$ (and e.g.~$a_i$ may be written $a[i]$) when they appear in a superscript or subscript.

\mn
8) $\bar x_{[u]} \defeq \LL x_i : i \in u\RR$

\mn
9) Forcing notions will be denoted by $\bbP$ and $\bbQ$. We adopt the Cohen convention that `$p \leq q$' means that $q$ gives more information (as {conditions in a forcing notion}).

\mn
10) $\unlhd$ means `is an initial segment,' and $\lhd$ means it is proper.

\mn
11) $\clT$ is a partial order or quasiorder, \underline{not} necessarily a tree. 

(Originally they were trees, but we later found it better to drop this --- see the end of \S2A. But it would be no problem to resurrect it in the future.)

\mn
12)  We use $\bfp$ to denote \emph{twinship parameters} (see Definition \ref{b3}) and $\bfm$ for {forcing examples}: see \S2B. 
\end{notation}

\mn
\begin{notation}\label{z6}
1) $\tau$ will denote a vocabulary: that is, a set of predicates and function symbols of finite \emph{arity} (that is, a finite number of places). Functions and individual constants are treated as predicates.

\mn
2) For models or structures, $\tau(M)$, $\tau(I)$, etc.~are defined naturally, as their vocabularies.

\mn
3) $\cL$ will denote a logic. $\bbL$ is first order logic, $\bbL_{\lambda,\mu}$ the usual infinitary logic.

$\cL(\tau)$ is the language: that is, a set of formulas $\varphi(\bar x)$ for the logic $\cL$ in the vocabulary $\tau$.

\mn
4) $T$ will denote a theory; complete first-order in the vocabulary $\tau_T = \tau(T)$, if not said otherwise. For simplicity, it will have elimination of quantifiers. (Particularly in \S0, we may forget to say `countable.')

\mn
5) For such $T$,
\begin{align*}
    \EC_\lambda(T) \defeq&\ \{M \models T : \|M\| = \lambda\}\\
    \EC(T) \defeq&\ \bigcup\limits_{\lambda \in \Card} \EC_\lambda(T).
\end{align*}
For $T_1 \supseteq T$,
\begin{align*}
    \PC_\lambda(T_1,T) \defeq&\ \{M \rest \tau_T : M \in \EC_\lambda(T_1)\}\\
    \PC(T_1,T) \defeq&\ \bigcup\limits_{\lambda \in \Card} \PC_\lambda(T_1,T).
\end{align*}
\end{notation}

\mn
\begin{definition}\label{z10}
For $\bbP$ a forcing notion, we define:
\begin{enumerate}
    \item $\kappa(\bbP) \defeq \min \{\kappa :\ \Vdash_\bbP\text{``there is a new } A \subseteq \kappa"\}$
\sn
    \item $\spec(\bbP) \defeq$ 
    \begin{align*}
        \big\{ (\kappa,\lambda,\clT) : &\ \clT \text{ is a subtree of ${}^{\kappa>}\!\lambda$ of cardinality $\lambda$ such that}\\
        &\text{ forcing with $\bbP$ adds a \underline{new} } \eta \in \lim_\kappa(\clT)  \big\}
    \end{align*}
    (I.e. $\eta \in {}^\kappa\lambda \setminus \bfV$ with $\eps < \kappa \Rightarrow \eta \rest \eps \in \clT$.)
\sn
    \item ``$(\kappa,\lambda) \in \spec(\bbP)$" will be shorthand for $\big(\exists \clT \big) \big[ (\kappa,\lambda,\clT) \in \spec(\bbP) \big]$.
\end{enumerate}
\end{definition}

\mn
\underline{\textbf{Question:}} will $\spec(\bbP)$ be interesting? Do we use
$\clT$ or just use a ``new" directed $\bfG \subseteq \bbP$?

\mn
\begin{convention}\label{z13}
If not stated otherwise, we assume $\bbP$ is such that
$$
\big(\forall p \in \bbP \big) \big[\kappa(\bbP_{\geq p}) = \kappa(\bbP) \big] 
$$
(where $\bbP_{\geq p} \defeq \bbP \rest \{q : q \geq_\bbP p\}$).
\end{convention}

\mn
This choice  can be justified by the following observation.
\begin{observation}\label{z14}

For any forcing notion $\bbP$ there is a maximal antichain $\bfI$ such that for all $p \in \bfI$, either $(a)$ or $(b)$ holds:
\begin{enumerate}[$(a)$]
    \item $\bbP_{\geq p}$ is a trivial forcing: i.e.~it is a directed quasiorder. (This means that any two members are \emph{compatible}, in that they have a common upper bound.)
\sn
    \item $\kappa(\bbP_{\geq p})$ is a well-defined regular uncountable cardinal $\leq \|\bbP_{\geq p}\| \leq \|\bbP\|$, and there is a minimal $\lambda \leq \|\bbP\|$ such that $(\kappa,\lambda) \in \spec(\bbP_{\geq p})$.
\end{enumerate}
\end{observation}

\begin{PROOF}{\ref{z14}}
Obvious, but we shall elaborate.

Let $\bfI$ be a maximal antichain of $\bbP$ included in 
\begin{align*}
    \{p \in \bbP : &\text{ either $\bbP_{\geq p}$ is directed \underline{or} there does not}\\
    &\text{ exist a $q \geq p$ such that $\bbP_{\geq q}$ is directed}\}.
\end{align*}
It will suffice to prove that this $\bfI$ is as promised. It will also suffice to deal with a non-directed $\bbP_{\geq p}$, where $p \in \bfI$.

\smallskip
Let $\mu \defeq \|\bbP_{\geq p}\|$ and let $f : \bbP_{\geq p} \to \mu$ be a bijection. Let $\name\eta \in {}^\mu2$ be the $\bbP$-name satisfying
$$
p \Vdash `\name\eta(\alpha) = 1\text{' \underline{iff} } p \in \Name\bfG_\bbP. 
$$
Let $\Lambda_1 \defeq \big\{ (\kappa,\lambda,\name\nu) :\ \Vdash_\bbP ``\name \nu \in {}^\kappa\lambda \text{ is not from $\bfV$, but } \eps < \kappa \Rightarrow \name\nu \rest \eps \in \bfV" \big\}$

Clearly $(\mu,\mu,\name\eta) \in \Lambda_1$. Hence there is a minimal $\kappa_*$ such that the set 
$$
\big\{ (\lambda,\name\nu) : (\kappa_*,\lambda,\name\nu) \in \Lambda_1\big\}
$$ 
is non-empty. Similarly, there is a minimal $\lambda_*$ such that $(\kappa_*,\lambda_*,\name\nu) \in \Lambda_1$ for some $\name\nu$. 
Let 
$$
\clT \defeq \big\{ \varrho \in {}^{\kappa_*>}\lambda_* : (\exists q \in \bbP_{\geq p}) \big[ q \Vdash ``\varrho \lhd \name\nu_*" \big] \big\}.
$$
Now $\clT,\kappa_*,\lambda_*,\name\nu$ are as promised in \ref{z10}.
\end{PROOF}

\mn
\begin{definition}\label{z16}
The following definition will be used mainly in \ref{c5}.
\begin{enumerate}[(a)]
    \item  $\tau(\mu,\kappa) = \tau_{\mu,\kappa}$ is the vocabulary with function symbols
\[
    \{F_{i,j} : i < \mu,\ j < \kappa\},
\]
    where $F_{i,j}$ is a $j$-place function symbol and $\kappa$ is a regular cardinal.
\sn
    \item  $\clM_{\mu,\kappa}(I)$ is the free $\tau_{\mu,\kappa}$-algebra generated by $I$.
\sn
    \item  We may write $\clM_\mu(I)$ when $\kappa = \aleph_0$, and $\clM(I)$ when $\mu = \kappa = \aleph_0$.
\end{enumerate}
\end{definition}

\mn
\begin{remark}\label{z29}
Concerning the first approach (see \S0B) we will define some games which witness the equivalence of two models in some strong logic.
\end{remark}

\mn
\begin{definition}\label{z32}
1) We say the models $M$ and $N$ are \emph{cofinally $(\lambda,\zeta)$-equivalent} \underline{when} there exist $\subseteq$-increasing sequences $\olsi M = \LL M_\alpha : \alpha < \lambda\RR$ and $\olsi N = \LL N_\alpha : \alpha < \lambda\RR$ satisfying the following.
\begin{enumerate}
    \item $M = \bigcup\limits_{\alpha<\lambda} M_\alpha$ and $N = \bigcup\limits_{\alpha<\lambda} N_\alpha$.
\sn
    \item The pro-isomorphism player \textsf{ISO} has a winning strategy in the game\\ $\Game_\zeta^\iso(\olsi M,\olsi N)$ defined below.
\end{enumerate}
A play of the game $\Game_\zeta^\iso(\olsi M,\olsi N)$ between the players \textsf{ISO} and \textsf{ANTI} lasts $\zeta$-many moves. In the $\eps^\tthh$ move, the \textsf{ANTI} player chooses $\alpha_\eps \in \big( \bigcup\limits_{\xi < \eps} \alpha_\xi,\lambda \big)$ and the \textsf{ISO} player responds with an isomorphism $f_\eps : M_{\alpha_\eps} \to N_{\alpha_\eps}$ extending $\bigcup\limits_{\xi < \eps} f_\xi$.

\mn
2) If $\|M\| = \|N\| = \lambda$ (and for transparency, both have universe $\lambda$) then we may demand that the isomorphism $f_\eps$ 
be a function from $\alpha_\eps$ onto $\alpha_\eps$.

\mn
3) In part (1), we may replace the ordinal $\zeta$ by a tree $\clT$ with $\lambda$-many levels and no $\lambda$-branch. 

By this we mean: in the $\eps^\tthh$ move of a play of 
$\Game_\clT^\iso(\olsi M,\olsi N)$, \textsf{ANTI} starts by choosing $t_\eps$, a member of the $\eps^\tthh$ level of $\clT$ which is $\lhd$-above $t_\xi$ for all $\xi < \eps$, and then $\alpha_\eps \in \big( \bigcup\limits_{\xi < \eps} \alpha_\xi,\lambda \big)$. Then \textsf{ISO} chooses $f_\eps : M_{\alpha_\eps} \to N_{\alpha_\eps}$ as before. A player loses if they have no legal move on their turn.

Note that $\alpha_\eps$ chosen exactly as in part (1), and does not depend on $t_\eps$. The tree simply functions as the game's `clock:' if \textsf{ISO} chooses a valid $f_\eps$ and \textsf{ANTI} has no valid $t_{\eps+1}$, then  \textsf{ISO} wins the play.
\end{definition}

\mn
E.g.~we have (in other variants of the game we get equivalence):
\begin{claim}\label{z35}
If $M$ and $N$ are two models of cardinality $\lambda \in \Reg$ and are cofinally $(\lambda,\omega)$-equivalent, \underline{then} they are $\bbL_{\infty,\lambda}$-equivalent.
\end{claim}

\medskip
The following property of linear orders will be used for proving that models are not isomorphic.
\begin{definition}\label{z38}
A model $J$ (usually a linear order) has the $\lambda$-\emph{indiscernibility} property when:
\begin{quotation}
     If $\bar t_\eps \in {}^{\omega>}\!J$ for $\eps < \lambda$, \underline{then} for some $A \in [\lambda]^\lambda$, the sequence $\LL \bar t_\eps : \eps \in A\RR$ is indiscernible for the quantifier-free formulas.
\end{quotation}
\end{definition}

\mn
\begin{fact}\label{z41}
If $\lambda$ is regular uncountable, \underline{then} any well-ordered set has the $\lambda$-indiscern-ibility property.

\end{fact}

\newpage

\section{GEM models}\label{S1}

Below, the reader may concentrate on $K_\oor$, $K_\org$, the order property, and the independence property.

\sn
Recall
\begin{notation}\label{a5}
1) Let $K$ denote a class of index models (i.e.~structures) which have the Ramsey property. (See \cite[1.10, p.330]{Sh:300}, \cite[1.15\subref{c2}]{Sh:E59}.) Members of $K$ will be denoted by $I$ and $J$; we shall use them for {constructing} generalized Ehrenfeucht-Mostowski models $\GEM(I,\Phi)$.   $\Phi$ (or $\Psi$) is called the \emph{blueprint}, and $\bfa = \LL \bar a_s : s \in I\RR$ will denote the \emph{skeleton}.

\mn
2) We may write $K_\x$ for (e.g.) $\x \in \{\oor,\org,\org(\bfn),\tr(\omega),\tr(\kappa),\tr(\bar\bfn),\oi(\partial)\}$.

In this case $K_\lambda^\x \defeq \{I \in K_\x : \|I\| = \lambda\}$.
\end{notation}

\medskip
Now we can define GEM models (\emph{Generalized Ehrenfeucht-Mostowski} models) for $K$. 
On this, see \cite[Ch.III, 1.6, p.329]{Sh:300} (revised in \cite[\S1B, 1.8\subref{b8}, p.9]{Sh:E59}). 

This usually requires generalizing Ramsey's theorem. Some examples of relevant classes:
\begin{example}\label{a8}
\begin{enumerate}
    \item $K = K_\oor$: the class of linear orders.
\sn
    \item $K = K_\omega^\tr = K_{\tr(\omega)}$: trees with $\omega+1$ levels. {We have} $P_i$ for $i \leq \omega$, $\lhd$ the tree order, and $<_\lex$ the lexicographical order. 
    
    (See \cite[1.9(4)\subref{b11}]{Sh:E59}.)
\sn
    \item [(B)$_\kappa$] $K = K_{\tr(\kappa)} = K_\kappa^\tr$: 
    similarly, but with $\kappa+1$ levels (so we have restriction functions 
    $\rest_{i,j}$). (See \cite[1.7(4), p.328]{Sh:300}, \cite[1.9(4)\subref{b11}]{Sh:E59}.)
\sn
    \item $K_\org$: linearly ordered graphs. 
    
    (See \ref{b19}, and more in \cite[1.18(5)\subref{c14}]{Sh:E59}.)
\sn
    \item $K_\dorg$ (\emph{directed ordered graphs}). Like $K_\org$, but the graph is directed. (We may also consider it as an {undirected} graph.)
\sn
    \item $K_{\org(\bfn)}$ for $\bfn \in [2,\omega]$ (see \cite{Sh:F2433}).
\sn
    \item $K_{\ptr(\sigma)} = K_\ptr^\sigma$ and $K_{\tr(\bar\bfn)} = K_{\tr(\bar\bfn)}^\sigma$ (for $\bar\bfn = \LL n_i : i < \sigma\RR \in {}^\sigma\omega$) are as in \cite[1.1(1)\subref{1.1}, 1.2\subref{1.2}]{Sh:511}, respectively.
    
    (Also called $K_{\str(\bar\bfn)}$; see more in \cite[Def. 5.1\subref{s1}, p.30]{Sh:F918}.)
\sn
    \item $K_{\oi(\gamma)}$; see \cite[Def. 2.1\subref{2b.1}]{Sh:897}.
\sn
    \item For any $K_\x$ such that $E,F_{\eta,\iota} \notin \tau(K_\x)$ for $\eta \in \clT_\bfp$ and $\iota = \pm1$, we define $K_{\clT,\iota}^\x$ --- the main case in this paper --- in \ref{b19}(1),(2).
\end{enumerate}
\end{example}

\mn
\begin{definition}\label{a41}
1) For $T$ a theory (not necessarily first-order), $K$ as in \ref{a5}(1), and $\kappa$ a cardinal, we define $\Upsilon_K[T,\kappa] = \Upsilon[T,\kappa,K]$ as the class of $K$-\emph{GEM blueprints} $\Phi$ (see \cite[Ch.III, 1.6, p.329]{Sh:300}, \cite[1.8\subref{b8}]{Sh:E59}).
\begin{enumerate}
    \item [$\boxplus$] For $I \in K$, $M = M_I \in \GEM(I,\Phi)$ 
    is a $\tau_\Phi$-model with skeleton $\bfa$.

    Pedantically, $(M,\bfa) \in \GEM(I,\Phi)$ satisfies
    \begin{enumerate}
        \item $\tau_\Phi$ is a vocabulary, and $M$ is a 
        $\tau_\Phi$-model.
\sn
        \item $\tau_\Phi \supseteq \tau_T$ and $\GEM_{\tau_T}(I,\Phi) = \GEM(I,\Phi) \rest \tau_T \models T$.
 \sn
        \item $\tau_\Phi$ is of cardinality $\leq \kappa$.
 \sn
        \item $\Phi$ is a function with domain 
        $$
        \QFT_K \defeq \big\{ \tp(\bar s,\varnothing,J) : J \in K,\ \bar s \in {}^{\omega>}\!J\big\},
        $$
        and $\Phi\big( \tp_\qf(\bar s,\varnothing,J)\big)$ is a complete $\bbL(\tau_\Phi)$-quantifier-free type.

        (`$\qf$' means \emph{quantifier-free}.)
\sn
        \item $\bfa = \LL \bar a_s : s \in I\RR$ is the so-called \emph{skeleton} of $M$.
 \sn
        \item $\lh(\bar a_s) = k_\Phi$ (So members of the skeleton are $k$-tuples. For simplicity, we will usually have $k_\Phi \defeq 1$. If so, we may write $a_s$ instead of $\LL a_s\RR$.)
 \sn
        \item $M$ is the closure of $\{\bar a_s : s \in I\}$.
 \sn
        \item $\bfa$ is \emph{qf-indiscernible} in $\GEM(I,\Phi)$, where `qf-indiscernible' is defined as in clause $\oplus_\qf$ below.
\sn
        \item If $\bar s \in {}^\eps\! I$ then $\ol{\bar a}_{\bar s} = \LL\ldots \caret \bar a_{s_\zeta} \caret\ldots\RR_{\zeta<\eps}$.

        (So if $k_\Phi \defeq 1$ then $\bar a_{\bar s} = \LL a_{s_\zeta} : \zeta < \eps\RR$.)
    \end{enumerate}
\mn
    \item [$\oplus_\qf$] \underline{If} $(s_0,\ldots,s_{n-1}), 
    (t_0,\ldots,t_{n-1}) \in {}^n\!I$ realize the same quantifier-free type in $I$, \underline{then} 
    $\bar a_{s_0} \caret \ldots \caret \bar a_{s_{n-1}}$ and 
    $\bar a_{t_0} \caret \ldots \caret \bar a_{t_{n-1}}$  realize the same quantifier-free type in $\GEM(I,\Phi)$. 
\end{enumerate}

\mn
1A) Of course, we are really interested in 
$\GEM_{\tau_T}(I,\Phi) = \GEM(I,\Phi) \rest \tau_T$.

\mn
1B) As implied above, we define $\GEM_\tau(I,\Phi) \defeq \GEM(I,\Phi) \rest \tau$ for $\tau \subseteq \tau_\Phi$.

\mn
2) We may write $\Upsilon_\kappa^K(T)$ of $\Upsilon_K(T,\kappa)$, and we may write $\Upsilon_\kappa^\x(T)$ for $K \defeq K_\x$ (e.g.~$\Upsilon_\kappa^\oor(T)$, $\Upsilon_\kappa^\org(T)$, 
$\Upsilon_\kappa^{\tr(\omega)}(T)$ for $K \defeq K_\oor,K_\org,K_{\tr(\omega)}$, respectively).
\end{definition}

\mn
The following definition also applies to non-first-order $T$ (and/or $\varphi$, or replace $\EC(T)$ by a class of models). When both are first-order, by the compactness theorem it suffices to use $\mu \defeq \aleph_0$.
\begin{definition}\label{a44}
1)
\begin{enumerate}
    \item We say that $\varphi = \varphi(\bar x_{[k]},\bar y_{[k]})$ \emph{witnesses that $T$ has the $\lepref{\lambda}$-order property} (not necessarily first-order) \underline{when} for every {$\mu < \lambda$}, there is $M \in \EC(T)$ and $\LL \bar a_\alpha : \alpha < \mu\RR \subseteq {}^\mu({}^k\!M)$ such that $M \models \varphi[\bar a_\alpha,\bar a_\beta]^{\iif(\alpha < \beta)}$.
\sn
    \item Let $T$ be first-order complete. We say $\varphi$ \emph{witnesses that $T$ is unstable} \underline{if} $\varphi$ is first-order and $T$ has the $\aleph_0$-order property as witnessed by $\varphi$.
\end{enumerate}

\mn
2)
\begin{enumerate}
    \item We say that $\varphi = \varphi(\bar x_{[k]},\bar y_{[k]})$ \emph{witnesses that $T$ has the $\lepref{\lambda}$-independence property} (not necessarily first-order) \underline{when} for every {$\mu < \lambda$} and graph $G = (\mu,R)$ on $\mu$, there are $M \in \EC(T)$ and $\LL \bar a_\alpha : \alpha < \mu\RR \subseteq {}^k\!M$ such that
    $$
    \alpha < \beta < \mu \Rightarrow \big[ M \models \varphi[\bar a_\alpha,\bar a_\beta] \Leftrightarrow G \models ``\alpha\ R\ \beta" \big].
    $$
    \item Let $T$ be first-order complete. We say $T$ is \emph{independent} (or \emph{has the independence property}) \underline{when} some first-order $\varphi(\bar x,\bar y) \in \bbL(\tau_T)$ witnesses the $\aleph_0$-independence property.

\end{enumerate}
\end{definition}

\mn
\begin{remark}\label{a50}
For all the examples in {\ref{a8}}, for the relevant (first-order) $T$, $\varphi$, or $\olsi\varphi$, there exists a suitable $\Phi$. (Usually we need that $K$ satisfies the Ramsey property; see \cite[III]{Sh:300} or \cite[\S1]{Sh:E59}.)
\end{remark}

E.g.

\begin{claim}\label{a51}
$1)$ If $T$ is first-order unstable (as witnessed by $\varphi = \varphi(\bar x,\bar y)$) and $T_1 \supseteq T$, \underline{then} there is $\Phi \in \Upsilon_\oor(T_1)$ such that: 
\begin{enumerate}[$(a)$]
    \item $|\tau_\Phi| = |T_1|$ 
\sn
    \item $k_\Phi = \lh(\bar x) = \lh(\bar y)$
 \sn
    \item $\GEM(I,\Phi) \models \varphi[\bar a_s,\bar a_t]^{\iif(s <_I t)}$ 
\end{enumerate}

\mn
$2)$ Let $\mu \geq \beth((2^\lambda)^+)$. If $T_1 \subseteq \bbL_{\lambda^+,\aleph_0}(\tau_1)$ is of cardinality $\leq\lambda$ and has the $\lepref{\mu}$-order property, as witnessed by $\varphi(\bar x_{[k]},\bar y_{[k]})$, \underline{then} the conclusion in part $(1)$ holds. (See more in \cite{Sh:E59}.)
\end{claim}





\mn
Note that the definition below formalizes the statements in \ref{a50}-\ref{a51}.
\begin{definition}\label{a56}
1) For $\Phi \in \Upsilon_\kappa[K]$, we say $\Phi$ \emph{represents} $(\varphi,R)$ \underline{when}:
\begin{enumerate}
    \item $R \in \tau(K)$ has arity $n$.
\sn
    \item $\varphi = \varphi(\bar x_0,\ldots,\bar x_{n-1}) \in \cL(\tau_\Phi)$ for some logic $\cL$, with $\lh(\bar x_\ell) = k_\Phi$.
\sn
    \item If $I \in K$, $(M,\bar\bfa) = \GEM(I,\Phi)$, and $\bar t \in {}^n\!I$, \underline{then}
    $$
    M \models \varphi[\bar a_{t_0},\ldots,\bar a_{t_{n-1}} ] \Leftrightarrow (t_0,\ldots,t_{n-1}) \in R^I.
    $$
\end{enumerate}


\mn
2) We may write $\varphi$ instead of $(\varphi,R)$ when $R$ is clear from the context. (E.g.~it is `$x<y$' for $K_\oor$, and $x\ R\ y$ in $K_\org$.)

\mn
3) Similarly {for} ``$\Phi$ represents $(\olsi\varphi,\olsi R)$," where $\olsi\varphi = \LL \varphi_\eps(\bar x_0,\ldots,\bar x_{m_\eps-1}) : \eps < \kappa\RR$.
\end{definition}

\medskip
The following definitions are natural ways to make demands even stronger than ``$M_1$ and $M_2$ are not isomorphic."
\begin{definition}\label{a59}
1) For $\tau$-models $M_1,M_2$ and $\varphi = \varphi( x_0,\ldots, x_{n-1}) \in \cL(\tau)$, we say $M_1$ is $(\lambda,\varphi)$-\emph{far from} $M_2$ \underline{when} there is a \emph{witness} $\LL a_\alpha : \alpha < \lambda\RR \in {}^\lambda (M_1)$. 

By \emph{this}, we mean: 
\begin{quotation}
    If 
    $\clU \in [\lambda]^\lambda$ and $(\forall \alpha \in \clU)[b_\alpha \in M_2]$, \underline{then} for some $\alpha_0 < \ldots < \alpha_{n-1}$ from $\clU$, we have
$$
    M_1 \models \varphi[a_{\alpha_0},\ldots,a_{\alpha_{n-1}}] \Leftrightarrow 
    M_2 \models \neg\varphi[ b_{\alpha_0},\ldots, b_{\alpha_{n-1}}].
$$
\end{quotation}

\mn
2) If $\varphi = \varphi(\bar x_0,\ldots,\bar x_{n-1}) \in \cL(\tau)$ with $(\forall \ell < n)[\lh(\bar x_\ell) = k]$ for some $k < \omega$, then above {we will write} $\bar a_\alpha \in {}^k\!M_1$, $\bar b_\alpha \in {}^k\!M_2$.

\mn
3) We say $M_1$ and $M_2$ \emph{are} $(\lambda,\varphi)$-\emph{far} \underline{when} $M_1$ is $(\lambda,\varphi)$-far from $M_2$ and $M_2$ is $(\lambda,\varphi)$-far from $M_1$.

\mn
4) For a logic $\cL$, $\tau$-models $M_1,M_2$, and $\varphi = \varphi(\bar x) \in \cL(\tau)$, we say $f$ is a $\varphi$-\emph{embedding of $M_1$ into $M_2$} \underline{when} for every $\bar a \in {}^{\lh(\bar x)}\!M_1$ we have
$$
M_1 \models \varphi[\bar a] \Leftrightarrow M_2 \models \varphi[f(\bar a)].
$$

\mn
5) ``$\olsi\varphi$-far" is defined similarly, for $\olsi\varphi = \LL \varphi_n : n \in [2,\bfn]\RR$.
\end{definition}

\mn
\begin{remark}\label{a65}
1) In \ref{a59}(1)-(5), if $M_1$ and $M_2$ are $\bbP$-twins (see Definition \ref{x1}(1), \ref{x2.0}(1)) and 
$$
\|\bbP\| < \lambda = \cf(\lambda) \leq \max(\|M_1\|, \|M_2\|),
$$ 
\underline{then} $M_1$ and $M_2$ cannot be `far' in any of those definitions.

\mn
[Why? Let $\bar a_\alpha \in {}^k\!M_1$ for $\alpha < \lambda$; let 
$$
\Vdash_\bbP``\name f : M_1 \to M_2 \text{ is an isomorphism (or just a $\varphi$-embedding)".}
$$ 
For each $\alpha < \lambda$ we can choose $p_\alpha \in \bbP$ and $\bar b_\alpha \in {}^k\!M_2$ such that $p_\alpha \Vdash ``\name f(\bar a_\alpha) = \bar b_\alpha$". As $\|\bbP\| < \lambda$, there is some $p_* \in \bbP$ such that the set $\{\alpha < \lambda : p_\alpha = p_*\}$ has cardinality $\lambda$. So recalling Definition \ref{a59}(1), $\LL \bar a_\alpha : \alpha < \lambda\RR$ cannot be a witness for ``$M_1$ is $(\lambda,\varphi)$-far from $M_2$."]

\medskip
So Definitions \ref{a68}(1),(2) below are an attempt to formulate notions of being `far' for which we might try to build examples.

\mn
2) In clause \ref{a68}(3)(C), we may demand that $\alpha_i$ is increasing {with $i$}.

\mn
3) Definition \ref{a68} will be used in \ref{c17} and its relatives.

\mn
4) Sometimes negation is not so handy. Then we may replace $\varphi$ by a pair $\LL \varphi_\Pos,\varphi_\Neg\RR$ in \ref{a59}+\ref{a68}, where $\varphi_\Pos = \varphi_\Pos(\bar x_0,\ldots,\bar x_{n-1})$ (and similarly for $\varphi_\Neg$).

What this means is that, in all occurrences,
\begin{enumerate}[$\bullet_1$]
    \item If $M_1 \models \varphi_\Pos[\bar a_0,\ldots]$ then $M_2 \models \varphi_\Pos[f(\bar a_0),\ldots]$, for $\bar a_0,\ldots,\bar a_{n-1} \in {}^k\!M_1$.
\sn
    \item If $M_1 \models \varphi_\Neg[\bar a_0,\ldots]$ then $M_2 \models \varphi_\Neg[f(\bar a_0),\ldots]$, for $\bar a_0,\ldots,\bar a_{n-1} \in {}^k\!M_1$.
\sn
    \item For no $\bar a_0^\ell,\ldots \bar a_{n-1}^\ell \in {}^k\!M_\ell$ (with $\ell \in \{1,2\}$) do we have 
    $$
    M_\ell \models \varphi_\Pos[\bar a_0^\ell,\ldots,\bar a_{n-1}^\ell] \wedge \varphi_\Neg[\bar a_0^\ell,\ldots,\bar a_{n-1}^\ell].
    $$
\end{enumerate}
\end{remark}

\mn
\begin{definition}\label{a68}
1) For $\tau$-models $M_1,M_2$ and $\varphi = \varphi( x_0,\ldots, x_{n-1}) \in \cL(\tau)$, we say $M_1$ is $(\lambda,\sigma,\varphi)$-\emph{far from} $M_2$ \underline{when} there is a {witness} $\LL \bfa_i : i < \sigma\RR$, where $\bfa_i = \LL a_{i,\alpha} : \alpha < \lambda\RR \in {}^\lambda (M_1)$. 

By \emph{this}, we mean: 
\begin{quotation}
    If $\clU_i \in [\lambda]^\lambda$ for $i < \sigma$, \underline{then} there does not exist a function 
    $$
    f : \{ a_{i,\alpha} : i < \sigma,\ \alpha \in \clU_i\} \to M_1
    $$
    preserving the satisfaction of $\varphi$.
\end{quotation}

\mn
2) Similarly for $\varphi = \varphi(\bar x_0,\ldots,\bar x_{n-1}) \in \cL(\tau)$ with $\lh(\bar x_0) = \ldots = \lh(\bar x_{n-1}) = k$.

\mn
3) We may replace $\varphi$ by $\Delta \subseteq \cL(\tau)$, or by $\LL \Delta_u : u \in [\sigma]^{< \aleph_0}\RR$, with the natural meaning:
\begin{enumerate}
    \item $\Delta_u \subseteq \big\{\varphi(\ldots,\bar x_i,\ldots)_{i \in u} : \varphi \in \cL(\tau) \big\}$
\sn
    \item $\lh(\bar a_{i,\alpha}) = \lh(\bar x_i) = k_i$
\sn
    \item There are $\clU_i \in [\lambda]^\lambda$ and $\bar b_{i,\alpha} \in {}^{k_i}(M_2)$ (for $i < \sigma$ and $\alpha \in \clU_i$) such that \underline{if} $u \in [\sigma]^{< \aleph_0}$, $\varphi(\ldots,\bar x_i,\ldots)_{i \in u} \in \Delta_u$, and $\alpha_i \in \clU_i$ for $i \in u$, \underline{then} 
    $$
    M_1 \models \varphi[\ldots,\bar a_{i,\alpha_i},\ldots]_{i\in u} \Leftrightarrow 
    M_2 \models \varphi[\ldots,\bar b_{i,\alpha_i},\ldots]_{i\in u}.
    $$
\end{enumerate}

\mn
4) We may add $\LL n_i : i < \sigma\RR$ with $n_i \leq \omega$, so $\varphi \in \Delta_u$ may have $\le n_i$ copies of $\bar x_i$ in our block of {arguments}.
\end{definition}

\newpage
\section{Toward $\bbP$-twins}\label{S2}

The idea below is that $\clT_\bfp$ is a forcing notion. However, sometimes we do not use the forcing notion we are interested in, but rather a derived one (e.g.~for Sacks forcing, we may use $\clT_\bfp \defeq ({}^{\omega>}2,\lhd)$). Even if $\clT_\bfp$ is {the forcing notion which interests us}, $\clB_\bfp$ (see \ref{b3}) will not necessarily be the family of all dense open subsets of $\clT_\bfp$.

\bn
\subsection{The Frame}

\begin{convention}\label{b2}
If not stated otherwise, $\bfp$ is a fixed weak twinship parameter (as in Definition \ref{b3}), although we will sometimes have different definitions for the tree-like and non-tree-like version. 

Earlier we thought it necessary to assume $\clT$ is a tree; now it seems this is not necessary, but neither option would be harmful to us so far.
\end{convention}

\mn
\begin{definition}\label{b3}
1) We say $\bfp$ is a \emph{weak twinship parameter} (or simply `a twinship parameter') \underline{when} it consists of:\footnote{
    So $\clT = \clT_\bfp$, etc. Hence we may write $\bfp = (\clT,\clB,\theta)$.
}
\begin{enumerate}
    \item A partial order $\clT$ such that any two elements $\eta,\nu \in \clT$ have a maximal lower bound\footnote{
        Many of the usual forcing notions fail this, but it will not be a problem to fix this. 
    } 
    (call it $\eta \wedge\nu$). 
\sn
    \item $\theta = \cf(\theta) \geq \aleph_0$

\sn
    \item $\clB$, a family of subsets of $\clT$ satisfying the following.
    \begin{enumerate}
        \item $\clB$ is closed under finite intersections, and each $D \in \clB$ is dense in $\clT$.
\sn
        \item If $\eta \in \clT$ then $\{\nu \in \clT : \eta <_\clT \nu \text{ or } \eta \perp \nu\} \in \clB$, or at least contains a member of $\clB$.
\sn
        \item If $\theta > \aleph_0$ \underline{then} the intersection of $<\theta$ many members of $\clB$ will always contain some other member of $\clB$.
    \end{enumerate}
\end{enumerate}

\mn
1A)  We say $\bfp$ is a \emph{strong} twinship parameter \underline{when} we add the following to clause (C):
\begin{enumerate}
    \item [(d)] No directed $\bfG \subseteq \clT$ meets every $D \in \clB$. 
\end{enumerate}

\mn
2) We say $\bfp$ is \emph{tree-like} \underline{when} in addition,
\begin{enumerate}
    \item $\clT$ is a tree with $\theta$-many levels.
\sn
    \item $(\forall \eta \in \clT)(\forall \eps < \theta)(\exists \nu \in \clT)[ \eta \leq_\clT \nu \wedge \lev_\clT(\nu) \geq \eps].$
\sn
    \item If $\eta \in D \in \clB$ \underline{then} $\eta < _\clT \nu \Rightarrow \nu \in D$. 
\end{enumerate}

\mn
3) We say $\bfp$ is \emph{well-founded} \underline{when} $\clT_\bfp$ has no infinite decreasing sequence.

\mn
4) So in part (1) we may write $\bfp = (\clT,\clB,\theta) = (\clT_\bfp,\clB_\bfp,\theta_\bfp)$, and in part (3) we may write $\bfp = (\clT,\clB) = (\clT_\bfp,\clB_\bfp)$.
\end{definition}

\mn
\begin{remark}\label{b7}
0) We would like to add the following demand to \ref{b3}(2), but it would hinder \ref{b71}:
\begin{enumerate}
    \item [(D)] $(\forall \eps < \theta)(\exists D \in \clB)(\forall \eta \in D)[\lev(\eta) \geq \eps]$.
\end{enumerate}

\mn
1) Originally the main cases were $K_\oor$ (the class of linear orders), $K_{\tr(\omega)}$ (trees with $\omega+1$ levels), and most importantly $K_\org$ (ordered graphs), but now it appears $K_{\tr(\omega)}$ is of limited importance.

\mn
2) Our aim, for a given $\bfp$, is to find non-isomorphic models (preferably \emph{far}, as well) for (e.g.)~a suitable complete first-order theory $T$ such that \underline{if} a forcing adds a new directed $\bfG \subseteq \clT$ meeting every $D \in \clB_\bfp$ (and if $\clT$ is a tree, then a new $\theta_\bfp$-branch of $\clT$ meeting every $D \in \clB_\bfp$)
\underline{then} they become isomorphic. 

Of course, if the forcing collapses some cardinal then this is trivially true (e.g.~for first-order countable $T$). We can rectify this by restricting ourselves to `interesting' $\bbP$-s (e.g.~ccc, $\Cohen_{\aleph_1}$ assuming $\CH$, or $\lepref{\theta}$-complete $\theta^+$-cc for a suitable $\theta = \theta^{<\theta}$) or requiring the models to be equivalent in some suitable logic. At any rate, to get something non-trivial, $T$ has to be somewhat complicated. This leads us to classification theory, so this is part of the classification program.
\end{remark}

\mn
\begin{definition}\label{b9}
1) We say $\bfG \subseteq \clT_\bfp$ \emph{solves} $\bfp$ when it is a directed subset meeting every $D \in \clB_\bfp$. 

\mn
2) We say $M$ and $N$ are $\bfp$-\emph{isomorphic} \underline{when} for {any} forcing notion $\bbP$, 
$$
\Vdash_\bbP ``\text{if some $\bfG \subseteq \clT_\bfp$ solves $\bfp$, then } M\cong N".
$$

\mn
3) We say $\bfp$ \emph{represents} the forcing notion $\bbP$ \underline{when} in any forcing extension $\bfV'$ of $\bfV$, $\bfp$ is solved in $\bfV'$ \underline{iff} in $\bfV'$ there exists a subset of $\bbP$ generic over $\bfV$.

\mn
4) We say `$M$ and $N$ are $\bfp$-twins' \underline{when} they are $\bfp$-isomorphic but not isomorphic.

\mn
5) For $\bfp$ a weak twinship parameter, we say the models $M$ and $N$ are \emph{strictly} $\bfp$-isomorphic \underline{when} for any forcing notion $\bbP$ we have
$$
\Vdash_\bbP ``\text{if some downward closed $\bfG \subseteq \clT_\bfp$ solves $\bfp$ and $\bfG \notin \bfV$, then } M\cong N".
$$
\end{definition}









\mn
\begin{definition}\label{b20}
1) We define
$$
\Omega_\bfp \defeq \big\{ \bfo = \big\LL (\eta_\ell,\iota_\ell) : \ell < k\big\RR : k < \omega,\ \eta_\ell \in \clT_\bfp,\ \iota_\ell = \pm 1\big\}.
$$
We may also denote $\bfo$ by the pair $(\bar\eta,\bar\iota) = (\bar\eta_\bfo,\bar\iota_\bfo) = \big(\LL \eta_\ell : \ell < k\RR, \LL \iota_\ell : \ell < k\RR\big)$.

As always, $\eta_\ell = \eta_{\bfo,\ell}$, $\iota_\ell = \iota_{\bfo,\ell}$, and $k = k_\bfo \defeq \lh(\bfo)$.

\mn
2) For $\bfo \in \Omega_\bfp$ and $\{F_{\eta_\ell,\iota_\ell} : \ell < k_\bfo\}$ a family of partial permutations of some set $\clU$, satisfying $F_{\eta,-\iota} = F_{\eta,\iota}^{-1}$, we define $F_\bfo$ naturally by induction on 
$\lh(\bfo)$:
\begin{itemize}
    \item $F_{\LL\ \RR} \defeq \id_\clU$ (I.e.~$F_{\LL\ \RR}(a) = a$ {iff} $a \in \clU$.)
\sn
    \item If $\bfo_2 = \bfo_1 \caret \big\LL (\eta,\iota)\big\RR$ then $F_{\bfo_2}(a) \defeq F_{\bfo_1}(F_{\eta,\iota}(a))$ (when it is well-defined).
\end{itemize}

\mn
3) For $\bfo$ and $F_{\eta_\ell,\iota_\ell}$ as above, the $\bfo$-\emph{orbit} of $a \in \clU$ is the sequence $\bar a = \LL a_\ell : \ell \leq k\RR$, where $a_k \defeq a$ and $a_\ell \defeq F_{\eta_\ell,\iota_\ell}(a_{\ell+1})$.

\mn
(Note that this sequence may not exist, as we do not require that 
$\dom(F_{\eta_\ell,\iota_\ell}) = \clU$.)

\mn
4) We say that the orbit $\bar a$ is \emph{reduced} \underline{when} 
$a_\ell \neq a_m$ for all $\ell < m \leq \lh(\bar a)$.

\mn
5) 
We say that $\bfo = \big\LL (\eta_\ell,\iota_\ell) : \ell < \lh(\bfo)\big\RR \in \Omega_\bfp$ is \emph{formally reduced} \underline{when}
\begin{itemize}
    \item If $\ell+1 < \lh(\bfo)$ and $\eta_\ell = \eta_{\ell+1}$, \underline{then} $\iota_\ell = \iota_{\ell+1}$.
\end{itemize}

\mn
6) $\Omega_\fr = \Omega_\bfp^\fr$ will denote the set of $\bfo \in \Omega_\bfp$ which are formally reduced.

\mn
7) $\Omega_\clS$ and $\Omega_\clS^\fr$ are defined similarly, replacing $\clT_\bfp$ by any set $\clS$.
\end{definition}

\mn
We first give a special case of the main definition, combining $K_\org$ with $\bfp$.
\begin{definition}[\textbf{Main Definition}]\label{b19}
0) $\tau_\org \defeq \{<,R\}$, where $<$ and $R$ are two-place {predicates}, and 
$$
K_\org \defeq \big\{ I =  \big(|I|,<^I,R^I  \big) : {<^I} \text{ is a linear order and $ \big(|I|,R^I \big)$ is a graph} \big\}.
$$

\mn
1) Let $\tau_\bfp^\org \defeq \{<,E,R,F_{\eta,\iota} : \eta \in \clT_\bfp,\ \iota = \pm 1\}$, where 
\begin{itemize}
    \item $<$ and $E$ are two-place {predicates}.
\sn
    \item $F_{\eta,\iota}$ is a unary function symbol (interpreted as a partial function).
\end{itemize}

\mn
1A) However, ``$x \in \dom(F_{\eta,\iota})$" is considered an atomic formula, and even\\ ``$x \in \dom(F_\bfo)$" for $\bfo \in \Omega_\bfp$.

\mn
2) Let $K_{\clT,0}^\org$ be the class of $\tau_\bfp^\org$-structures $I$ such that:
\begin{enumerate}
    \item $<_I$ is a linear order on $I$.
\sn
    \item For $\eta \in \clT$, $F_{\eta,\iota} = F_{\eta,\iota}^I$ is a partial automorphism of $(|I|,<_I,R^I)$, increasing with $\eta$, satisfying the following:
    \begin{enumerate}
        \item $F_{\eta,-\iota} = F_{\eta,\iota}^{-1}$.
\sn
        \item If $a \in I$ and $\iota = \pm1$, then $D_{\iota,a}^I \defeq \big\{ \eta : \dom(F_{\eta,\iota}) \ni a\big\} \in \clB$.

        (If $\iota=1$ we may omit it.)
\sn
        \item $F_\bfo^I$ is well-defined and a partial automorphism for all $\bfo \in \Omega_\bfp$.
\sn
        \item If $\eta \leq_\clT \nu$ and $\iota = \pm 1$ then $F_{\eta,\iota}^I \subseteq F_{\nu,\iota}^I$.
\sn
        \item \textbf{[Follows:]} If $\iota = \pm1$, $F_{\eta_\ell,\iota}^I(s_\ell) = t_\ell$ for $\ell = 1,2$, and $\eta_1,\eta_2$ are $\leq_\clT$-compatible, \underline{then}
$$
        s_1\, R^I\, s_2 \Leftrightarrow t_1\, R^I\, t_2.
$$
    \end{enumerate}
\sn
    \item $(|I|,R^I)$ is a graph.
\sn
    \item $E^I$ is the closure of 
    $$
    \big\{ \big(a,F_{\eta,\iota}^I(a)\big) : a \in \dom(F_{\eta,\iota}^I),\ \eta \in \clT_\bfp,\ \iota = \pm1\big\}
    $$
    to an equivalence relation.
\end{enumerate}

\mn
3) Let $K_{\clT,1}^\org$ be the class of $I \in K_{\clT,0}^\org$ such that
\underline{if} $\bfo \in \Omega_\fr$, $\lh(\bfo) = k \geq 1$, $a_k \in I$, and $a_0,\ldots,a_k$ is an $\bfo$-orbit, \underline{then} $a_k \neq a_0$.

\mn
4) For $I,J \in K_{\clT,0}^\org$, let `$I \subseteq J$' mean ${<_I} = {<_J} \rest |I|$, $R^I = R^J \rest |I|$, and $F_\bfo^I = F_\bfo^J \rest I$ for all $\bfo \in \Omega_\bfp$. 
\end{definition}

\mn
Similarly, we can define
\begin{definition}\label{b22}
Let $K = K_\x$ be as in \ref{a5}(1A) and $E,F_{\eta,\iota} \notin \tau(K_\x)$, where\\ $\{F_{\eta,\iota} : \eta \in \clT,\ \iota = \pm1\}$ are unary function symbols and $E$ a binary predicate. (Below, if $\clT = \clT_\bfp$ then we may replace $\clT$ by $\bfp$.)

\mn
1) Let $\tau_\clT^\x \defeq \tau(K_\x) \cup \{F_{\eta,\iota} : \eta \in \clT,\ \iota = \pm 1\} \cup \{E\}$.

\mn
2) $K_{\clT,0}^\x$ is the class of structures $J$ such that:
\begin{enumerate}
    \item 
    \begin{enumerate}
        \item $J$ is a $\tau_\clT^\x$-model.
\sn
        \item $J \rest \tau(K_\x) \in K_\x$
\sn
        \item Every $F_{\eta,\iota}^J$ is a partial automorphism of $J \rest \tau(K_\x) \in K_\x$, increasing with $\eta \in \clT$.
    \end{enumerate} 
\sn
    \item Clauses \ref{b19}(2)(B)(a)-(d) all hold, but (e) is replaced by:
    \begin{enumerate}
        \item [(e)$'$] {For} $R$ any predicate from $\tau(K_\x)$,
        if $F_{\eta_\ell,\iota_\ell}^I(s_\ell) = t_\ell$ for $\ell < \arity(R)$ and $\{\eta_\ell : \ell < \arity(R)\}$ has a common upper bound, \underline{then} 
$$
        \big\LL s_\ell : \ell < \arity(R) \big\RR \in R^I \Rightarrow \big\LL t_\ell : \ell < \arity(R) \big\RR \in R^I.
$$
    \end{enumerate}
    (Recall that we are assuming $\tau(K_\x)$ has only predicates:  function symbols and individual constants will be treated as predicates.)
\sn
    \item [(C),(D)] As in Definition \ref{b19}(2)(C),(D).
\end{enumerate}

\mn
3-4) As in \ref{b19}(3),(4).
\end{definition}

\mn
\begin{definition}\label{b34}
1) Let $\leq_\Omega$ be the following two-place relation on $\Omega_\bfp$.
\begin{enumerate}
    \item [$\circledast$] $\bfo_1 \leq_\Omega \bfo_2$ \underline{iff} ($\bfo_1,\bfo_2 \in \Omega_\bfp$ and)
    \begin{enumerate}
        \item $\lh(\bfo_1) = \lh(\bfo_2)$
\sn
        \item $\bar\iota^1 = \bar\iota^2$ (where $\bfo_\ell = (\bar\eta^\ell,\bar\iota^\ell)$).
\sn
        \item $\eta_k^1 \unlhd \eta_k^2$ for all $k < \lh(\bfo_1) = \lh(\bfo_2)$.
    \end{enumerate}
\end{enumerate}

\mn
2) $\bfo_1 \parallel_\Omega \bfo_2$ will be shorthand for ``$\bfo_1$ and 
$\bfo_2$ are compatible" (i.e.~have a common $\leq_\Omega$-upper bound). 
\end{definition}

\mn
\begin{definition}\label{b37}
1) For $J \in K_{\clT,1}^\x$ and $s \in J$, let $\Omega_s^J \defeq \{\bfo \in \Omega_\bfp^\fr : \dom(F_\bfo^J) \ni s\}$.

\mn
2) For $D \in \clB$ (see \ref{b3}(1)(C)), 
let 
$$
\Omega_D \defeq \big\{\bfo \in \Omega_\bfp^\fr : \{\eta_{\bfo,\ell} : \ell < \lh(\bfo)\} \subseteq D \big\}.
$$

\mn
3) Let $K_{\clT,2}^\x = K_{\bfp,2}^\x$ be the class of $J \in K_{\clT,1}^\x$ such that 
\begin{enumerate}
    \item $(\forall s \in J)(\exists D \in \clB)[ \Omega_s^J = \Omega_D]$
\sn
    \item If $F_\bfo^J(s) = t$ \underline{then} there exists $\bfo' \leq_\Omega \bfo$ such that $F_{\bfo'}^J(s) = t$ and
    $$
    (\forall \bfo'' <_\Omega \bfo') \big[s \notin \dom(F_{\bfo''}^J) \big].
    $$
\end{enumerate}
\end{definition}

\mn
Note that
\begin{fact}\label{b41}
1) If $\bfp$ (i.e.~$\clT_\bfp$) is well-founded, then clause (B) of \ref{b37}(3) always holds.

\mn
2)  If $\bfp$ it tree-like \underline{then} it is well-founded.
\end{fact}

\begin{PROOF}{\ref{b41}}
1) This follows immediately from $\bfp$ being well-founded --- 
$\{\bfo' : \bfo' \leq_\Omega \bfo\}$ has no infinite decreasing sequence for any $\bfo \in \Omega$.

\mn
2) Easy as well.
\end{PROOF}

\mn
\begin{remark}
The assumption that $\bfp$ is well-founded is not a serious hindrance, thanks to \ref{b71} below.
\end{remark}

\mn
\begin{claim}\label{b44}
$1)$ $\leq_\Omega$ is a partial order on $\Omega_\bfp$.

\mn
$2)$ $\bfo_1 \parallel_\Omega \bfo_2$ \underline{iff} $\lh(\bfo_1) = \lh(\bfo_2)$ and for $\ell < \lh(\bfo_1)$, $\iota_{\bfo_1,\ell} = \iota_{\bfo_2,\ell}$ and $\eta_{\bfo_1,\ell}$, $\eta_{\bfo_2,\ell}$ are $\leq_\clT$-compatible.

\mn
$3)$ Assume $\bfp$ is tree-like.

Any $\bfo_1,\bfo_2 \in \Omega_\bfp$ have a maximal common $\leq_\Omega$-lower bound; we will denote it by $\bfo_1 \wedge \bfo_2$.
\end{claim}

\begin{PROOF}{\ref{b44}}
Easy. (E.g.~part (4) follows from \ref{b3}(A).)
\end{PROOF}

\mn
\begin{claim}\label{b47}
$1)$ If $J \in K_{\clT,0}^\x$ and $\bfV_{\!1} \defeq \bfV^\bbP$ (or an extension), \underline{then} $(A) \Rightarrow (B)$, where
\begin{enumerate}[$(A)$]
    \item $\bfG \in \bfV_{\!1}$ solves $\bfp$: that is, it is a directed subset of $\clT_\bfp$  such that 
$$
    (\forall D \in \clB_\bfp)[\bfG \cap D \neq \varnothing].
$$
    \item $F_\bfG \defeq \bigcup\limits_{\eta \in \bfG} F_{\eta,1}^J$ is an automorphism of $J \rest \tau_\org$. (So in the proof of \emph{\ref{c17}(1)}, it will induce an isomorphism from one distinguished subset of $J \rest \tau_\org$ to another, and hence from $M_1^+$ to $M_2^+$.)
\end{enumerate}

\mn
$2)$ If $\bfp$ is tree-like, then clause $(A)$ is equivalent to
\begin{enumerate}
    \item [$(A)'$]
    \begin{enumerate}[$\bullet_1$]
            \item  $\clT_\bfp \subseteq {}^{\theta>}\!\lambda$ is a tree, $\bfV_{\!1} \models ``\eta \in {}^\theta\!\lambda"$, and $(\forall \eta < \theta)[\eta \rest \eps \in \clT]$.
\sn
            \item $\bfG \defeq \{\eta \rest \eps : \eps < \theta\}$ solves $\bfp$. 
        \end{enumerate}
\end{enumerate}
\end{claim}

\begin{PROOF}{\ref{b47}}
First,
\begin{enumerate}
    \item [$(*)_1$] $F_\bfG$ is a well-defined function.
\end{enumerate}
[Why? As $\bfG \subseteq \clT_\bfp$, each $F_\eta^J$ is a function. Furthermore, if $\eta ,\nu \in \bfG$ \underline{then} there is $\rho \in \bfG$ such that $\eta \leq_\clT \rho \wedge \nu \leq_\clT \rho$ (because $\bfG$ is directed) and $F_\eta^J \subseteq F_\rho^J \wedge F_\nu^J \subseteq F_\rho^J$ (recalling $J \in K_\bfp^\x$ and \ref{b19}(2)(B)(c)).] 
\begin{enumerate}
    \item [$(*)_2$] $F_\bfG$ is a partial automorphism of $J \rest \tau_\x$.
\end{enumerate}
[Why? Similarly to the proof of $(*)_1$.]
\begin{enumerate}
    \item [$(*)_3$] $\dom(F_\bfG) = J$.
\end{enumerate}
[Why? Recall \ref{b19}(2)(B)(b), and that $\bfG$ meets every $D \in \clB_\bfp$.]

\mn
Lastly,
\begin{enumerate}
    \item [$(*)_4$] $\rang(F_\bfG) = J$.
\end{enumerate}
[Why? Like $(*)_3$, recalling $\dom(F_{\eta,1}^J) = \rang(F_{\eta,-1}^J)$.]

Together we are done.
\end{PROOF}

\mn
\begin{observation}\label{b50}
For every $D \in \clB_\bfp$ there exists $I \in K_{\clT,2}^\org$ such that:
\begin{enumerate}[$(a)$]
    \item $\Omega_s^I = \Omega_D$ for every $s \in I$.
\sn
    \item For every $s \in I$, $I = \{F_\bfo^I(s) : \bfo \in \Omega_D\}$. 
\end{enumerate}
\end{observation}

\begin{PROOF}{\ref{b50}}
Straightforward.
\end{PROOF}

\vspace{8mm}
\subsection{Examples}
First we present examples of $\bfp$ with $\theta_\bfp = \aleph_0$ (so they do not fit the theorems in \S3).

\begin{claim}[\textbf{Example / Claim}]\label{b53}\ 

\mn
$1)$ The following example is a strong tree-like twinship parameter witnessing the Cohen Real forcing.

Let $\bfp = \bfp_\Cohen = \bfp[\Cohen]$ consist of:\footnote{
    So $\clT = \clT_\bfp$, etc.
}
\begin{enumerate}[$(a)$]
    \item $\clT_\bfp = ({}^{\omega>}2,\lhd)$ (so it is tree-like, and $\theta = \aleph_0$).
\sn
    \item $\clB_\bfp$, the set of open dense subsets of ${}^{\omega>}2$.
\end{enumerate}

\mn
$2)$ If $\lambda \defeq \cov(\meagre)$, then for some $\clB \subseteq \clB_{\bfp[\Cohen]}$ of cardinality $\lambda$, the pair 
$\big(({}^{\omega>}2,\lhd),\clB\big)$ is a strong tree-like twinship parameter --- that is, it satisfies Definition \emph{\ref{b3}}.

\mn
$3)$ Let $\clT$ be $({}^{\omega>}2,\lhd)$, \emph{or} a subtree of $({}^{\delta>}2,\lhd)$ for some $\delta$ such that 
$$
(\forall \eta \in \clT)(\forall \eps < \delta)( \exists \nu \in \clT \cap {}^\eps2)[ \nu \unlhd \eta \vee \eta \unlhd \nu].
$$
Define $\bfp = \bfp[\clT]$ by $\clT_\bfp \defeq \clT$ and $\clB_\bfp \defeq \{ \clT \setminus {}^{\eps\geq}2 : \eps < \delta \}$. 

\underline{Then} $\bfp$ is a weak twinship parameter.

\mn
$4)$ All of the examples above are well-founded.
\end{claim}

\begin{PROOF}{\ref{b53}}
1) Covered in the proof of part (2).

\mn
2) Let $\LL Z_\alpha : \alpha < \lambda\RR$ be a sequence of meagre subsets of ${}^\omega2$ such that $\bigcup\limits_{\alpha < \lambda}Z_\alpha = {}^\omega2$. Let $Z_\alpha \subseteq \bigcup\limits_n Z_{\alpha,n}$, where $Z_{\alpha,n}$ is a closed nowhere dense subset of ${}^\omega2$, $\subseteq$-increasing with $n$. Let $\clB \defeq \{D_{u,n} :u \in [\lambda]^{< \aleph_0},\ n < \omega \}$, where
$$
D_{u,n} \defeq \Big\{\eta \in {}^{\omega>}2 : \neg \big(\exists \rho \in \textstyle\bigcup\limits_{\alpha \in u} Z_{\alpha,n} \big) \big[ \eta \lhd \rho \big] \Big\}.
$$

This obviously suffices, but let us elaborate, checking clauses (A)-(C) of Definition \ref{b3}. 

\mn
\textbf{Clause (A):} As $\clT \defeq ({}^{\omega>}2,\lhd)$, this is clear.

\mn
\textbf{Clause (B):} $\theta \defeq \aleph_0$, so it is regular and infinite.

\mn
\textbf{Clause (C):} 

\begin{enumerate}
    \item [$(*)_1$] For $\alpha < \lambda$ and $n < \omega$, $\clT_{\alpha,\{n\}} \defeq \{ \eta \rest n : \eta \in Z_{\alpha,n} \}$
    is a nowhere dense subtree of ${}^{\omega>}2$.\\ (That is, $\big( \forall \eta \in {}^{\omega>}2 \big) \big( \exists \nu \in {}^{\omega>}2 \setminus \clT_{\alpha,\{n\}} \big) \big[ \eta \lhd \nu \big]$.)
\end{enumerate}
[Why? Because $Z_{\alpha,n}$ is a nowhere dense subset of ${}^\omega2$.]

\sn
\begin{enumerate}
    \item [$(*)_2$] Each $D_{u,n}$ is a nowhere dense subtree of ${}^{\omega>}2$.
\end{enumerate}
[Why? By $(*)_1$.]

\sn
\begin{enumerate}
    \item [$(*)_3$] $\clB$ is a family of subsets of $\clT \defeq {}^{\omega>}2$.
\end{enumerate}

So to verify clause (C) we have to check the sub-clauses.

\mn
\textbf{(C)(a):} $\clB$ is closed under finite intersections.

Let $m < \omega$ and $\{D_\ell : \ell < m\} \subseteq \clB$. So for each $\ell < m$ there are $n_\ell < \omega$ and $u_\ell \in [\lambda]^{\aleph_0}$ such that $D_\ell = D_{u_\ell,n_\ell}$. 

Let $n \defeq \max(\{n_\ell : \ell < m\} \cup \{1\})$ and $u \defeq \bigcup\limits_{\ell < m} u_\ell$. Easily, $D_{u,n} \in \clB$ and $D_{u,n} = \bigcap\limits_{\ell < m} D_{u_\ell,n} \subseteq \bigcap\limits_{\ell < m}  D_{u_\ell,n_\ell} = \bigcap\limits_{\ell < m} D_\ell$ are as required.

\mn
\textbf{(C)(b):} Let $m < \omega$ and $\eta \in \clT$. So for every $\nu \in {}^n2$ we can find $\rho_\nu \in {}^\omega2$ such that $\nu \lhd \rho_\nu$. By the choice of $\LL Z_\alpha : \alpha < \lambda\RR$, there exists an $\alpha_\nu < \lambda$ such that $\rho_\nu \in Z_{\alpha_\nu}$. Hence $\rho_\nu \in Z_{\alpha_\nu,m_\nu}$ for some $m_\nu < \omega$.

Let $m \defeq \max\{m_\nu : \nu \in {}^n2\}$, so 
$$
\nu \in {}^n2 \Rightarrow \rho_\nu \in Z_{\alpha_\nu,m_\nu} \subseteq Z_{\alpha_\nu,m} \subseteq \bigcap\limits_{\nu' \in {}^n2} Z_{\alpha_{\nu'},m} = D_{u,m},
$$
where $u \defeq \{\alpha_{\nu'} : \nu' \in {}^n2\}$.

Obviously $D_{u,m} \in \clB$ and is disjoint to ${}^{n\geq}2$, so it is as required.

\mn
\textbf{(C)(c):} Vacuous.

\mn
\textbf{Clause (C)(d):} (of \ref{b3}(1A)) 

If $\bfG$ solves $\bfp$, then there exists $\eta \in {}^\omega2$ such that 
$$
(\forall \alpha < \lambda)(\forall m < \omega)(\forall^\infty n)[ \eta \rest n \notin D_{\{\alpha\},m}],
$$
but this contradicts our choice of $\LL Z_\alpha : \alpha < \lambda\RR$.

\mn
3-4) Easy.
\end{PROOF}

\mn
\begin{claim}[\textbf{Example / Claim}]\label{b56}
The following example is a strong well-founded tree-like twinship parameter which witnesses the Random Real forcing.

Let $\bfT_*$ be the set of sequences $\olsi \clT = \LL \clT_n : n < \omega\RR$ such that:
\begin{enumerate}[$(A)$]
    \item $\clT_n$ is a subtree of ${}^{\omega>}2$.
\sn
    \item $\Leb(\bigcup\limits_{n \in [m,\omega)}\lim\clT_n) = 1$ for every $m < \omega$.
\sn
    \item Let $\eta_n \defeq \tr(\clT_n)$ (the \emph{trunk} of $\clT_n$). The sequence $\LL \eta_n : n < \omega\RR$ is without repetition, and 
    $$
    \eta_n \lhd \eta_m \notin \clT_n \Rightarrow (\exists \rho \in {}^{\omega>}2 \setminus \clT_n)[\eta_n \lhd \rho \lhd \eta_m].
    $$
\end{enumerate}
Let $\bfT \defeq \big\{\LL \clT_{\alpha,n} : n < \omega\RR : \alpha < \lambda \big\} \subseteq \bfT_*$ be such that 
$$
\bigcap\limits_{\alpha < \lambda}\bigcup\limits_{n < \omega} \lim \clT_{\alpha,n} = \varnothing
$$
(Clearly there is such a sequence of length $\lambda$ when $\lambda \defeq \cov(\nll)$.)

Let $\bfp = \bfp_\Random = \bfp_{\Random(\bfT)}$ consist of:
\begin{enumerate}[$(a)$]
    \item $\clT_\bfp = ({}^{\omega>}2,\lhd)$ (so $\theta = \aleph_0$).
\sn
    \item $\clB_\bfp \defeq \{D_{u,n} : u \in [\lambda]^{< \aleph_0},\ n < \omega \}$, where
$$
    D_{u,n} \defeq \Big\{\eta \in {}^{\omega>}2 : \lh(\eta) \geq n,\ (\forall \alpha \in u)(\exists m)\big[\tr(\clT_{\alpha,m}) \unlhd \eta \in \clT_{\alpha,m}\big] \Big\}
$$
\end{enumerate}
\end{claim}

\begin{PROOF}{\ref{b56}}
Similar to the proof of \ref{b53}.
\end{PROOF}

\medskip
The following wide family of examples cover Random Real forcing, Cohen forcing, and virtually every forcing which adds a new set of ordinals (even e.g.~Prikry forcing).

\begin{definition}\label{b59}
1) We say $\bfm = (\lambda,\theta,\clT,\bbP,\name\eta) = (\lambda_\bfm,\theta_\bfm,\ldots)$ is a \emph{forcing example} \underline{when}:
\begin{enumerate}
    \item $\lambda \geq 2$ and $\theta = \cf(\theta) \geq \aleph_0$.
\sn
    \item $\clT$ is a subtree of $({}^{\theta>}\!\lambda,\lhd)$ (so it is closed under initial segments).
\sn
    \item $\bbP$ is a forcing notion which preserves ``$\theta$ is regular."
\sn
    \item $\name\eta$ is a $\bbP$-name of a member of ${}^{\theta}\!\lambda$.
\sn
    \item $\Vdash_\bbP ``\name\eta \notin \bfV \text{ and } (\forall\eps < \theta)[\name\eta_\eps \in \clT]"$, where $\name\eta_\eps \defeq \name\eta \rest \eps \in {}^\eps\!\lambda$.
\sn
    \item  For transparency, we demand $(\forall \nu \in \clT)(\exists p \in \bbP)[p \Vdash``\nu \lhd \name\eta"]$.
\end{enumerate}

\mn
2) For $\bfm$ as above, let $\bfp = \bfp_\bfm = (\clT_\bfm,\clB_\bfm,\theta_\bfm)$ be defined as follows.
\begin{enumerate}
    \item $\clT_\bfp \defeq \clT_\bfm$, $\theta_\bfp \defeq \theta_\bfm$.
\sn
    \item $\clB_\bfm \defeq \big\{D_{\bfI,f} \subseteq \clT : (\bfI,f) \in \set(\bbP)\big\}$, \underline{where}
    \begin{enumerate}
        \item $\set(\bbP)$ is the set of pairs $(\bfI,f)$ such that:
        \begin{enumerate}
            \item $\bfI$ is a maximal antichain of $\bbP$, or its completion.
\sn
            \item $f : \bfI \to \clT$
\sn
            \item $f(p) = \nu \Rightarrow p \Vdash_\bbP ``\nu \lhd \name\eta$"
        \end{enumerate}
\sn
        \item $D_{\bfI,f} \defeq \big\{ \nu \in \clT : (\exists p \in \bfI)
        [f(p) \unlhd \nu]
        \big\}$.
    \end{enumerate}
\end{enumerate}

\mn
3) For a forcing notion $\bbP$, we define $\bfp[\bbP]$ as $\clT_\bfp \defeq \bbP$, $\clB_\bfp$ the family of dense open subsets of $\bbP$, and $\theta_\bfp \defeq \{\theta : \text{forcing with $\bbP$ adds a new function } f \in {}^\theta\Ord\}$.
\end{definition}

\mn
\begin{remark}\label{b60}
In \ref{b59}(1), we can use $\bfp_\bfm' = (\clT_\bfm,\clB_\bfm',\theta_\bfm)$ (so $\clT_\bfm \subseteq {}^{\theta>}\!\lambda$), where 
$\clB_\bfm'$ is a subset of $\clB_\bfm$ {containing} $\{ \clT \setminus {}^{\eps\geq}\lambda : \eps < \delta \}$, 
closed under finite intersections, with \underline{no} solution. (See \ref{b9}.)

Then there would be no $\eta \in {}^\theta\!\lambda$ such that 
$$
(\forall D \in \clB)(\forall^\theta \eps < \theta)[\eta \rest \eps \in D].
$$
\end{remark}

\mn
\begin{claim}\label{b62}
$0)$ If $\bbP$ is a forcing notion, and $\lambda$ and $\kappa$ are minimal such that 
$$
\Vdash_\bbP \text{\emph{``there is a new} }\eta \in {}^\kappa\!\lambda",
$$
\underline{then} there exists $\clT \subseteq {}^{\kappa>}\!\lambda$ of cardinality $\leq \|\bbP\|$ and a forcing example $\bfm$ such that
$$
(\lambda_\bfm,\theta_\bfm,\clT_\bfm,\bbP_\bfm) = (\lambda,\kappa,\clT,\bbP).
$$

\mn
$1)$ If $\bfm$ is a forcing example \underline{then} $\bfp_\bfm$ is a {tree-like} strong twinship parameter. 


\mn
$2)$ If $\bfm$ is a forcing example and $\bbP_\bfm$ is $\lepref{\aleph_1}$-complete (or at least adds no new $\omega$-sequence of ordinals), \underline{then} necessarily $\theta_\bfm > \aleph_0$.

\mn
$3)$ If $\bfm = (\lambda,\theta,\clT,\bbP,\name\eta)$ is a forcing example and $D \in \clB_{\bfp_\bfm}$, \underline{then}
\begin{enumerate}[$\bullet_1$]
    \item $\Vdash_\bbP (\forall^\infty \eps < \theta) [\name\eta \rest \eps \in D]$
\sn
    \item $\bfp_\bfm$ has a solution in $\bfV^\bbP$.
\end{enumerate}

\mn
$4)$ If $\bbP$ is a non-trivial forcing \underline{then} $\bfp[\bbP]$ (see \emph{\ref{b59}(3)}) is a strong twinship parameter and $\Vdash_\bbP ``\bfp[\bbP] \text{\emph{ has a solution"}}$.

\mn
$5)$ If $\clT$ is a Suslin tree, \underline{then} there exists a forcing example $\bfm$ with $\bbP_\bfm \defeq \clT$ and $\theta_\bfm \defeq \aleph_1$.

\mn
$6)$ All of the $\bfp$-s defined above are well-founded.
\end{claim}

\begin{PROOF}{\ref{b62}}
Easy.

E.g.~in part (1), to verify Definition \ref{b3}(1)(C)(c), use 
$$
\Vdash_\bbP ``\theta \in \Reg"
$$ 
from \ref{b59}(1)(C).
\end{PROOF}

\mn
\begin{remark}\label{b63}
In \ref{b62}(3)\,$\bullet_2$, we may use a smaller $\clB' \subseteq \clB_{\bfp[\bbP]}$ as in previous examples. (The gain from this is that if $\clB'$ has a smaller cardinality, then \ref{c11} and \ref{c14} apply to more values of $\lambda$.)


\end{remark}

\mn
\begin{claim}\label{b65}
$1)$ Assume $\bfm$ is a forcing example.
\begin{enumerate}[$(a)$]
    \item If $\bbP_\bfm$ is Cohen forcing then $\theta_\bfm = \aleph_0$.
\sn
    \item If $\bbP_\bfm$ is $\lepref{\aleph_1}$-complete then $\theta_\bfm \geq \aleph_1$.
\sn
    \item If $\bbP_\bfm$ adds no new $\omega$-sequence of ordinals then $\theta_\bfm \geq \aleph_1$.
\end{enumerate}

\mn
$2)$ Assume $\bbP$ is a non-trivial forcing (i.e.~above every $p \in \bbP$ there are two incompatible members). Let $\kappa(\bbP)$ be as in \emph{\ref{z10}}, \emph{\ref{z13}}.

\underline{Then} for some $\lambda \leq |\bbP|$, there is a forcing example 
$\bfm$ with $\bbP_\bfm = \bbP$, $\theta_\bfm = \kappa(\bbP)$, and $\lambda_\bfm = \lambda$.
\end{claim}

\begin{PROOF}{\ref{b65}}
Easy, but we elaborate.

\mn
1) \textbf{Clause $(a)$:} Follows from the definition of Cohen forcing, which is $({}^{\omega>}2,\lhd)$ or something equivalent.

\mn
\textbf{Clause $(b)$:} By \ref{b59}(1)(B),(D), because forcing by $\bbP$ adds no new $\omega$-sequence of members of $\bfV$.

\mn
\textbf{Clause $(c)$:} Similarly.

\mn
2) Follows from Observation \ref{z14}.
\end{PROOF}

\mn
\begin{claim}\label{b68}
Assume  $\bfm$ is a forcing example with $\theta_\bfm > \aleph_0$.

\mn
$1)$ $\bbP_\bfm \neq \mathsf{Sacks}$. In fact, Sacks forcing adds no solution to $\bfp$.

\mn
$2)$ If $\bbQ$ is a definition of a forcing notion, is non-trivial and ccc, and {the truth values of} $``p \leq_\bbQ q"$,  $``p \perp_\bbQ q"$, and $``\bbQ$ is ccc$"$ are preserved in $\bfV^\bbQ$, \underline{then} $\bbP_\bfm \neq \bbQ$. 
Moreover, $\bbQ$ adds no solution to $\bfp$.

\mn
\emph{2A)} If $\theta_\bfm > \aleph_1$ then the preservation assumption in part $(2)$ may be omitted.

\mn
$3)$ $\bbP_\bfm$ fails the $\theta_\bfm$-Knaster condition (e.g.~$\Random$).
\end{claim}

\sn
\begin{remark}\label{b69}
Similarly for the forcing notions from Ros\l anowski-Shelah \cite{Sh:470} (and \cite{Sh:949}) --- see more in \cite{Sh:F2433}..
\end{remark}

\sn
\begin{PROOF}{\ref{b68}}
Let $\bfp \defeq \bfp_\bfm$, $\bbP \defeq \bbP_\bfm$, etc.

Towards contradiction, suppose 
$$
\Vdash_\bbP ``\name\nu \in \lim (\clT_\bfp) \text{ is a new $\theta_\bfp$-branch".}
$$

\mn
1) 
Let $\name\nu \in {}^\omega2$ be the generic real for Sacks forcing, and $\name \eta \defeq \name\eta_\bfm$. Clearly for a dense set of $p \in \Sacks$, $p$ forces that for some $u = u_p \in [\theta]^{\aleph_0}$ we can compute $\name\nu$ from $\name\eta \rest u_p$. 

So for every $p$ in this dense set, as $u$ is uncountable and $\theta$ is a regular uncountable cardinal, clearly $\zeta_p \defeq \bigcup\limits_{\alpha \in u_p} \alpha+1 < \theta$. As $\name\eta_\bfm \rest n = (\name\eta \rest \zeta_p) \rest u$, clearly
we have $p \Vdash_\bbP``\name\eta \rest \zeta_p$ is a new sequence". This gives us our contradiction.

\bn
2) Toward contradiction, assume $p \in \bbQ$ and $p \Vdash_\bbQ ``\name\eta \text{ solves } \bfp_\bfm$" (so it is a $\theta_\bfm$-branch of 
$\clT_{\bfp_\bfm}$ meeting every $D \in \clB_{\bfp_\bfm}$).

For every $\eps < \theta$, let 
$$
\Lambda_\eps \defeq \{\nu \in \clT_{\bfp_\bfm} : \lh(\nu) = \eps, \text{ and some $q \in \bbQ$ above $p$ forces } ``\nu\lhd \name\eta"\}.
$$
As $\bbQ$ satisfies the countable chain condition, clearly $\Lambda_\eps$ is countable and 
\begin{enumerate}
    \item [$(*)_1$] $(\forall \eps < \zeta < \theta)[\nu \in \Lambda_\zeta \Rightarrow \nu \rest \eps \in \Lambda_\eps]$.
\end{enumerate}
Also, 
\begin{enumerate}
    \item [$(*)_2$] $(\forall \eps < \zeta < \theta_{\bfp_\bfm}) (\forall \nu \in \Lambda_\eps)(\exists \rho \in \Lambda_\zeta)[\nu \lhd \rho]$.
\end{enumerate}
Moreover, as $\name\eta$ solves $\bfp_\bfm$, we have
$\Vdash_\bbQ ``\name\eta \notin \bfV$", hence 
\begin{enumerate}
    \item [$(*)_3$] $\big(\forall \eps < \theta_{\bfp_\bfm}\big) \big( \forall \nu \in \Lambda_\eps \big) \big( \exists \zeta \in [\eps,\theta_{\bfp_\bfm}) \big) \big( \exists \rho_1,\rho_2  \in \Lambda_\zeta\big) \big[ \rho_1 \neq \rho_2  \wedge \nu \lhd \rho_1 \wedge \nu \lhd \rho_2\big].$
\end{enumerate}
(Recall that $\bbQ$ satisfies the ccc.) This implies $\Lambda = \bigcup\limits_{\eps < \theta_\bfp} \Lambda_\eps$ is a tree with no $\theta_\bfp$-branches.

Let $\bfG \subseteq \bbQ$ be generic over $\bfV$ and contain $p$. Now in $\bfV[\bfG]$ we have a $\theta_{\bfp}$-branch $\eta$ of $\Lambda$, and for every $\eps < \theta_\bfp$ there exist $\xi = \xi_\eps \in (\eps,\theta_\bfp)$ and ${\varrho_\eps} \in \Lambda_{\xi} \setminus \{\eta \rest \xi\}$. 
\begin{enumerate}
    \item [$(*)_4$] So for some unbounded $A \subseteq \theta_\bfp$ we have 
    $$
    (\forall \eps,\zeta \in A)[ \eps < \zeta \Rightarrow \xi_\eps < \zeta].
    $$
\end{enumerate}

\sn
For every $\eps \in A$, there exists $p_\eps \in \bfG$ such that
\begin{enumerate}
    \item [$(*)_5$] $\bfV \models \big[p_\eps \Vdash``\varrho_\eps \lhd \name\eta" \big]$.
\end{enumerate}

\sn
So for all $\eps \neq \zeta \in A$, 
\begin{enumerate}
    \item [$(*)_6$]
    \begin{enumerate}
        \item $\bfV \models ``\varrho_\eps \text{ and $\varrho_\zeta$ are $\lhd$-incomparable in }\clT_\bfm"$
\sn
        \item $\bfV \models ``p_\eps \text{ and $p_\zeta$ are incompatible in }\bbQ"$.
    \end{enumerate}
\end{enumerate}
[Why? Note that $\lh(\varrho_\eps) = \xi_\eps < \zeta \leq \lh(\varrho_\zeta)$, so if $\varrho_\eps$ and $\varrho_\zeta$ are $\lhd$-compatible (in $\clT_\bfp$) then necessarily $\varrho_\eps \lhd \varrho_\zeta$. But $\lh(\varrho_\eps) = \xi_\eps$ and $\varrho_\eps \neq \eta \rest \xi_\eps \lhd \varrho_\zeta$, so if $q \in \bbQ$ is above $p_\eps$ and $p_\zeta$, then we get a contradiction.]

\smallskip
Note that each $p_\eps$ belongs to $\bbQ^\bfV \subseteq \bfV$, but $\bar p_A = \LL p_\eps : \eps \in A\RR$ may not be in $\bfV$ (just in $\bfV[\bfG]$).
So $\bar p_A$ contradicts our assumption that forcing with $\bbQ$ preserves ``$\bbQ$ satisfies the ccc and $p \perp q$" (i.e.~$p$ and $q$ are incompatible).


\mn
2A) Using the choices in the proof of {part (2)}, consider 
$\Lambda = \bigcup\limits_{\eps<\theta} \Lambda_\eps$. So it is a subtree of ${}^{\theta>}\!\lambda$, and each level is countable and non-empty. It is known that no such tree exists. 

(Alternatively: clearly without loss of generality $\|\bbP\| \geq \theta$. Let $\bbQ \in \clH(\chi)$; let $N \prec (\clH(\chi),\in)$ be of cardinality $<\theta$, $\bar\Lambda = \LL \Lambda_\eps : \eps < \theta\RR \in N$, and demand that $\delta \defeq N \cap \theta$ have uncountable cofinality. Choose $\eta \in \Lambda_\delta$ and continue as in the proof of part (2).)

\mn
3) For $\eps < \theta_\bfm$, let $p_\eps \in \bbP_\bfm$ force a value to $\name\eta_\eps$ (call it $\nu_\eps$). If $\bbP_\bfm$ satisfies the $\theta_\bfm$-Knaster condition, then there exists a set $\clU \in [\theta_\bfm]^{\theta_\bfm}$ such that $\LL p_\eps : \eps \in \clU\RR$ are pairwise compatible. Hence 
$$
\eps < \zeta \in \clU \Rightarrow \nu_\eps \lhd \nu_\zeta
$$
and $\nu \defeq \bigcup\limits_{\eps \in \clU} \nu_\eps$ is a branch of $\clT_{\bfp_\bfm}$.

Replacing $p_\eps$ by $p_{\min(\clU \setminus \eps)}$ for all $\eps \in \theta_\bfm \setminus \clU$, without loss of generality $\nu_\eps = \nu \rest \eps$ for every $\eps < \kappa$. It suffices to prove that 
$$
(\forall D \in \clB_{\bfp_\bfm}) \big[ D \cap \{\nu_\eps : \eps < \kappa\} \neq \varnothing \big].
$$
For every $D \in \clB_{\bfp_\bfm}$ and each $\eps < \kappa$, there exists $q_\eps \in \bbP_\bfm$ above $p_\eps$ forcing $\varrho_\eps \lhd \name\eta$ for some $\varrho_\eps \in D$, and we let $\zeta_\eps \defeq \max\{\eps,\lh(\varrho_\eps)\}$.

Now there necessarily exist $\eps_1,\eps_2 < \kappa$ such that $\zeta_{\eps_1} < 
{\eps_2}$ and $q_{\eps_1},q_{\eps_2}$ have a common upper bound. Let $r$ be such an upper bound.
So 
$$
r \Vdash ``\varrho_{\eps_1},\nu_{\eps_2} \text{ are both} \unlhd\, \name\eta",
$$ 
hence $\varrho_{\eps_1}$ and $\nu_{\eps_2}$ are $\unlhd$-compatible. But as 
$\lh(\varrho_{\eps_1}) \leq \zeta_{\eps_1} \leq \eps_2$, necessarily 
$\varrho_{\eps_1} \unlhd \nu_{\eps_2}$. But $\varrho_{\eps_1} \in D \in \clB_{\bfp_\bfm}$, which is an open dense subset of $\clT_\bfp$, hence $\nu_{\eps_2} \in D$. Therefore $D \cap \{\nu_\eps : \eps < \kappa\} \neq \varnothing$.

As $D$ was an arbitrary member of $\clB_\bfp$, the set $\{\nu_\eps : \eps < \kappa\}$ solves $\bfp_\bfm$ --- a contradiction.
\end{PROOF}

\bn
The demand `$\bfp$ is well-founded or tree-like'\footnote{
    See \ref{b3}(2),(3).
} 
may seem unreasonable, but
\begin{claim}\label{b71}
Assume $\bfp$ is a $[$weak/strong$]$ twinship parameter.
\begin{enumerate}[$(a)$]
    \item For $\eta \in \clT_\bfp$, let 
    $$
    \delta(\eta,\bfp) \defeq \sup\!\big\{\alpha+1 : (\exists \bar\nu \in {}^\alpha\clT_\bfp)[\bar \nu \text{ is increasing}\; \wedge\, \nu_0 = \eta] \big\}.
    $$
    \item Let $\bfI$ be a maximal antichain of $\clT_\bfp$ such that 
    $$
    {(\forall \eta \in \bfI)}(\forall \nu \geq_\clT \eta) \big[ \delta(\eta,\bfp) = \delta(\nu,\bfp) \big].
    $$
    \item For $r \in \bfI$, let ${\bfp_r} = (\clT_r',\clB_r',\theta_r')$ be defined as follows.
    \begin{enumerate}[$\bullet_1$]
        \item $\clT_r$ is the set of $<_\clT$-increasing sequences of length 
        a successor ordinal.
\sn
        \item ${<_{\clT_r}} \defeq {\lhd}$ (`is an initial segment of . . .').
\sn
        \item $\clB_r \defeq \{D_{[r]} : D \in \clB\}$, where 
        $$
        D_{[r]} \defeq \{\bar p \in \clT_r : \bar p \cap D \neq \varnothing\}.
        $$
        \item $\theta_r \defeq \theta$.
    \end{enumerate}
\end{enumerate}
\underline{Then} for each $r \in \bfI$:
\begin{enumerate}
    \item[$(d)$] $\bfp_r$ is a $[$weak/strong$]$ twinship parameter.
\sn
    \item[$(e)$] $\bfp_r$ is well-founded, and even tree-like (except that the {set of levels} of $\clT_{\bfp_r}$ is $\delta(r,\bfp)$, which is not necessarily a cardinal).
\sn
    \item[$(f)$] $|\clB_{\bfp_r}| \leq |\clB_\bfp|$
\sn
    \item[$(g)$] For any forcing extension $\bfV^\bbQ$ of $\bfV$, $\bfp$ has a solution \underline{iff} there exists $r \in \bfI$ such that $\bfp_r$ has a solution.
\end{enumerate}
Note: If $\delta(-,\bfp)$ is constant and $0$ is the minimal element of $\clT_\bfp$, then we can always choose $\bfI \defeq \{0\}$.
\end{claim}

\begin{PROOF}{\ref{b71}}
Straightforward, using \ref{z14}.
\end{PROOF}


\bn
\subsection{How `nice' are the classes $K_\x$?}

\sn
(Note that \ref{b26}--\ref{b29} will not be used later.) 

\begin{observation}\label{b26}
$1)$ Each of the classes from \emph{\ref{a8}} is an AEC (see \emph{\ref{q5}}) with the JEP (``Joint Embedding Property") and amalgamation, see e.g. \cite[Def.~2.5]{Sh:88}. (For $K_{\tr(\kappa)}$, recall that the unique member of $P_0^I$ is an individual constant, so identified.)

\mn
$2)$ This also applies for $K_{\clT,\iota}^\oor$ and $K_{\clT,\iota}^\org$ (for $\iota = 0,1,2$).
\end{observation}

\begin{PROOF}{\ref{b26}}
Straightforward.
\end{PROOF}

\mn
\begin{definition}\label{b31}
1) Let $\bfS_\clT^\oor$ be the {class} of $I \in K_\clT^\oor$ such that for any $s \in I$, 
$$
\{F_\bfo^I(s) : \bfo \in \Omega_\bfp,\ F_\bfo^I(s) \text{ is well-defined}\}
$$ 
is equal to the set of elements of $I$.

\mn
2) For $\bfS \subseteq \bfS_\clT^\oor$, let $K_\clT^\oor[\bfS]$ be the class of $I \in K_\clT^\oor$ such that for every $s \in I$, 
$$
I \rest \big\{F_\bfo^I(s) : \bfo \in \Omega_\bfp,\ F_\bfo^I(s) \text{ is well-defined}\big\}
$$
is isomorphic to some member of $\bfS$.

\mn
3) $\bfS_\clT^\org$ and $K_\clT^\org[\bfS]$ are defined similarly.

\mn
4) For {any} $K$, we can define $K_\clT$, $\bfS_\clT^K$, and $K_\clT[\bfS]$.
\end{definition}

\mn
\begin{claim}\label{b29}
$K_\oor$, $K_\org$, $K_{\clT,\ell}^\oor$, and $K_{\clT,\ell}^\org$ are universal
 classes.   
That is, if $M$ is a $\tau(K_\bfS)$-model and every finitely generated submodel belongs to $K_\bfS$, then $M \in K_\bfS$.
\end{claim}

\begin{PROOF}{\ref{b29}}
Obvious.
\end{PROOF}

\mn
Quoting \cite[1.2\subref{a5}]{Sh:88r}:
\begin{definition}\label{q5}  
We say $\gk$ is a AEC with 
LST number $\lambda(\gk) = \LST_\gk$ \underline{if}: 

\mn
\textbf{Ax.0}:  The truth of `$M \in K$' and `$N \le_\gk M$' 
depends on $N$ and $M$ only up to isomorphism; i.e.
$$
M \in K \wedge M \cong N \Rightarrow N \in K
$$ 
and 
`if $N \le_\gk M$, $f$ is an isomorphism
from $M$ onto the $\tau$-model $M'$, and $f \rest N$ is an isomorphism
from $N$ onto $N'$, \underline{then} $N' \le_\gk M'$.'

\mn
\textbf{Ax.I}:  if $M \le_\gk N$ 
then $M \subseteq N$ (i.e.~$M$ is a submodel of $N$).

\mn
\textbf{Ax.II}:  $M_0 \le_\gk M_1 \le_\gk M_2$ 
implies $M_0 \le_\gk M_2$ and $M \le_\gk M$ for $M \in K$.

\mn
\textbf{Ax.III}:  If $\lambda$ is a regular cardinal, $M_i$
is $\le_\gk$-increasing (i.e.~$i < j < \lambda$ implies 
$M_i \le_\gk M_j$) and
continuous (i.e.~for $\delta < \lambda$, $M_\delta = \bigcup\limits_{i < \delta}
M_i$) for $i < \lambda$ \underline{then} 
$$M_0 \le_\gk \bigcup\limits_{i < \lambda} M_i.$$

\mn
\textbf{Ax.IV}:  If $\lambda$ is a regular cardinal and $M_i$ 
(for $i < \lambda)$ is $\le_\gk$-increasing continuous and
$M_i \le_\gk N$ for $i < \lambda$ \underline{then} 
$\bigcup\limits_{i < \lambda} M_i \le_\gk N$.

\mn
\textbf{Ax.V}:  If $N_0 \subseteq N_1 \le_\gk M$ and $N_0 \le_\gk M$ 
\underline{then} $N_0 \le_\gk N_1$.

\mn
\textbf{Ax.VI}:  If $A \subseteq N \in K$ and $|A| \le \LST_\gk$, 
then for some $M \le_\gk N$, we have $A \subseteq |M|$ and $\|M\| \le 
\LST_\gk$ (and $\LST_\gk$ is the minimal
infinite cardinal satisfying this axiom which is $\ge |\tau|$; 
the $\ge |\tau|$ is for notational simplicity). 
\end{definition}

\newpage
\section{On Existence For independent $T$}\label{S3}

\begin{convention}\label{c2}
$\bfp$ is a (weak) twinship parameter (that is, as in Definition \ref{b3}).
\end{convention}

\mn
\begin{definition}\label{c5}
Assume $\lambda > \kappa + \aleph_0$ and $\kappa \geq 2$ such that $\alpha < \lambda \Rightarrow |\alpha|^{< \kappa} < \lambda$.

\mn
1) We say a graph $G$ is $(\lambda,\kappa)$-\emph{entangled} \underline{when} we have (A) $\Rightarrow$ (B), where
\begin{enumerate}
    \item 
    \begin{enumerate}
        \item $\eps < \kappa$
\sn
        \item $\bar a_\alpha = \LL a_{\alpha,\zeta} : \zeta < \eps\RR \in {}^\eps G$, and each $\bar a_\alpha$ is without repetitions\\ 
        (for $\alpha < \lambda$).
\sn
        \item For all $\alpha \neq \beta < \lambda$, the sets 
        $\{a_{\alpha,\zeta} : \zeta < \eps\}$ and 
        $\{a_{\beta,\zeta} : \zeta < \eps\}$ are disjoint.
    \end{enumerate}
\sn
    \item For every $X \subseteq \eps \times \eps$, there exist $\alpha < \beta < \lambda$ such that
    $$
    \big(\forall \zeta,\xi < \eps \big) \big[ a_{\alpha,\zeta}\, R^G\, a_{\beta,\xi}\ \Leftrightarrow\ (\zeta,\xi) \in X\big].
    $$
\end{enumerate}

\mn
2) We say $I \in K_{\clT,0}^\org$ (of cardinality $\geq \lambda$) is $(\lambda,\kappa)$-{entangled} \underline{when} we have\\ (A) $\Rightarrow$ (B), where:
\begin{enumerate}
    \item As above, but adding:
    \begin{enumerate}
        \item [(d)] If $\zeta,\xi < \eps$, $\bfo \in \Omega_\bfp$, and $\alpha < \beta < \lambda$, \underline{then} 
        $$
        F_\bfo^I(a_{\alpha,\zeta}) = a_{\alpha,\xi} \Leftrightarrow F_\bfo^I(a_{\beta,\zeta}) = a_{\beta,\xi}.
        $$

        \item [(e)] If $\alpha < \beta < \lambda$ and 
        $\gamma < \delta < \lambda$, \underline{then}
        $$
        (\forall\zeta,\xi < \eps)[a_{\alpha,\zeta} <_I a_{\beta,\xi} \Leftrightarrow a_{\gamma,\zeta} <_I a_{\delta,\xi}].
        $$
    \end{enumerate}
\sn
    \item For every $X \subseteq \eps \times \eps$, there exist $\alpha < \beta < \lambda$ such that
    $$
    \big(\forall \zeta,\xi < \eps \big) \big[ a_{\alpha,\zeta}\, R^I\, a_{\beta,\xi}\ \Leftrightarrow\ (\zeta,\xi) \in X\big],
    $$
    \emph{provided that}
    \begin{itemize}
        \item If $\gamma < \lambda$, $\bfo_\ell \in \Omega$, $F_{\bfo_\ell}^I(a_{\gamma,\zeta_\ell}) = a_{\gamma,\xi_\ell}$ for $\ell = 1,2$, {and $\bfo_1,\bfo_2$ have a common $\leq_\Omega$-upper bound}, then 
$$
        (\zeta_1,\zeta_2) \in X \Leftrightarrow (\xi_1,\xi_2) \in X.
$$
    \end{itemize}
\end{enumerate}


\mn
3) If we omit $\kappa$ and simply write `$\lambda$-entangled,' we mean $\kappa \defeq \aleph_0$.
\end{definition}


\medskip
We will state the following for a cardinal $\kappa$, but as in \cite{Sh:E59} $\kappa \defeq \aleph_0$ if not stated otherwise.
\begin{definition}\label{c8}
Assume $K$ is as in \ref{a5}, $\tau_{\mu,\kappa}$ is as in \ref{z16} (quoting from \cite{Sh:E59}), and $\Sigma$ is a set of $\tau_{\mu,\kappa}$-terms 
$\sigma(\bar x)$, where $\bar x = \LL x_\zeta : \zeta < \eps\RR$ for some 
$\eps < \kappa$. Further assume that $|\Sigma| \leq \mu$. (If $\mu = \mu^{< \kappa}$, this is automatic. If in addition $\Sigma$ is the set of all $\tau_{\mu,\kappa}$-terms, then we may omit $\Sigma$.)

{Then} for $I,J \in K$, we say that $I$ is \emph{strictly} $(\mu,\kappa)$-$\Gamma$-$\Sigma$-\emph{unembeddable} into $J$ (we may write `$\mu$' instead of $(\mu,\aleph_0)$)
\underline{when} we have `(A) $\Rightarrow$ (B),' where:
\begin{enumerate}
    \item 
    \begin{enumerate}
        \item $\clM_{\mu,\kappa}(J)$ is a $\tau_{\mu,\kappa}$-structure as in \ref{z16}.
\sn
        \item $F : I \to \clM_{\mu,\kappa}(J)$
\sn
        \item For each $s \in I$, $F(s)$ is of the form $\sigma_s(\bar t_s)$, where
        \begin{enumerate}
            \item $\sigma_s \in \Sigma$ (which is a set of $\tau_{\mu,\kappa}$-terms).
\sn
            \item 
            $\lh(\bar t_s) = \eps(\sigma_s) = \eps_s < \kappa$
\sn
            \item $\bar t_s = \LL t_{s,\eps} : \eps < \eps_s\RR \in {}^{\eps_s}\!J$.
\sn
            \item If $K$ has a linear order and 
        $$
            \kappa > \aleph_0 \Rightarrow J \text{ is well-ordered,}
        $$ 
            then $\bar t_s$ is $<_J$-increasing.
        \end{enumerate}
\sn
        \item $\Gamma$ is a set of pairs $\big( p_1(\bar x),p_2(\bar x) \big)$ such that for some $\eps < \kappa$, $I' \in K$, and $\bar t_\ell \in {}^\eps(I')$ (for $\ell = 1,2$), we have 
        $$
        p_\ell \subseteq \tp_\qf(\bar t_\ell,\varnothing,I').
        $$
    \end{enumerate}
\sn
    \item There exist $\eps < \kappa$, $\bar s_1,\bar s_2 \in {}^\eps\!I$, and $(p_1,p_2) \in \Gamma$ such that:
    \begin{enumerate}
        \item $p_1 \subseteq \tp_\qf(\bar s_1,\varnothing,I)$ and $p_2 \subseteq \tp_\qf(\bar s_2,\varnothing,I)\big)$.
\sn
        \item $\sigma_{s_{1,\zeta}} = \sigma_{s_{2,\zeta}}$ for all $\zeta < \eps$.
\sn
        \item The sequences $(\bar t_{s_{1,0}} \caret \ldots \caret \bar t_{s_{1,\zeta}} \caret \ldots)_{\zeta < \eps}$ and $(\bar t_{ s_{2,0}} \caret \ldots \caret \bar t_{ s_{2,\zeta}} \caret \ldots)_{\zeta < \eps}$ realize the same quantifier-free types in $J$.
    \end{enumerate}
\end{enumerate}
\end{definition}


\bn
The following claims will be used to get twins for independent theories in \ref{c17}.

\sn
\begin{claim}\label{c11}
Suppose $\lambda$ is regular,  
$\mu \in [\aleph_0, \lambda)$, $\bfp$ is a strong twinship parameter, $\theta_\bfp > \aleph_0$, and $|\clT_\bfp|^+ + |\clB_\bfp| \leq \lambda$.
Let $\Sigma$ be as in \emph{\ref{c8}}.

\underline{Then} we have \emph{`$\boxplus_1 \Rightarrow\boxplus_2$'}, where:
\begin{enumerate}[$\boxplus_1$]
    \item 
    \begin{enumerate}[$(a)$]
        \item $J \in K_{\clT,2}^\org$ (see Definition \emph{\ref{b37}(3)}) is of cardinality $\lambda$. 
\sn
        \item $\Gamma_{\!\org} \defeq \big\{ \big( p_1(x_0,x_1), p_2(x_0,x_1)  \big)\big\}$, where 
        $$
        p_1(x_0,x_1) \defeq \big[ x_0 < x_1\ \wedge\ x_0\, R\, x_1 \big]
        $$
        and 
        $$
        p_2(x_0,x_1) \defeq \big[ x_0 < x_1\ \wedge\ \neg(x_0\, R\, x_1) \big].
        $$
        \item $X \in [J]^\lambda$ is well-ordered and maximal such that 
        $x \neq y \in X$ implies $y \notin \cl_J(\{x\})$
        (equivalently, $\neg[x\ E^J\, y]$).
\sn
        \item $I \defeq (J \rest \tau_\org) \rest X \in K_\org$ is 
        $\lambda$-entangled.
\sn        
        \item For $\ell = 1,2$, we define 
        $$
        X_\ell \defeq \big\{F_\bfo^I(a) : a \in X,\ \bfo \in \Omega_a^J,\text{ and } \lh(\bfo) = \ell \mod 2 \big\}.
        $$
        \item 
        \begin{enumerate}    
            \item If $D \in \clB_\bfp$ then $|Y_D^0| \geq \lambda$, where $Y_D^0 \defeq \{s \in X : \Omega_s^J = \Omega_D\}$ (recalling Definition \emph{\ref{b37}}).
\sn
            \item $X$ is the disjoint union of $\LL Y_D^0 : D \in \clB_\bfp\RR$.
\sn
            \item (Note $X = X_1 \cup X_2$.)
        \end{enumerate}
\sn
        \item $I_\ell \defeq I \rest X_\ell$ for $\ell = 1,2$.
\sn
        \item $R^J = \big\{ \big(F_\bfo^J(s),F_\bfo^J(t) \big) : \bfo \in \Omega,\  s \neq t \in X \cap \dom(F_\bfo^I),\ (s,t) \in R^I \big\}$
\sn
        \item For $s_1,s_2 \in X$, if $\bfo_\ell \in \Omega_{s_\ell}^J$ for $\ell = 1,2$, \underline{then}
        $$
        s_1 <_J s_2 \Leftrightarrow F_{\bfo_1}^J(s_1) <_J F_{\bfo_2}^J(s_2).
        $$
        \item If $t,v \in X$ and $\bfo_1,\bfo_2 \in \Omega_t^J  = \Omega_v^J$, then 
        $$
        F_{\bfo_1}^J(t) <_J F_{\bfo_2}^J(t) \Leftrightarrow F_{\bfo_1}^J(v) <_J F_{\bfo_2}^J(v).
        $$
    \end{enumerate}
\sn
    \item $I_1$ is strictly $\mu$-$\Gamma_{\!\org}$-$\Sigma$-{unembeddable} into $I_2$.
\end{enumerate}
\end{claim}

\begin{PROOF}{\ref{c11}}
This is a special case of \ref{c14} proved below, where we choose 
$Y \defeq X$ and $Z \defeq X_2$.
(So $I_2 \defeq I \rest X_2$ here is equal to $I \rest Z$ there, and $I_1 \defeq I \rest X_1$ here \emph{contains} $I \rest Y$ from there, hence the conclusion of \ref{c14} implies $\boxplus_2$ here.)

Now we have to verify that the conditions of \ref{c14} hold. 
This is straightforward, noting that in clause (d)\,$\bullet_2$, the 
$\lambda$-indiscernibility property follows from $J$ being well-ordered, by \ref{c11}$\boxplus_1$(a) (recalling \ref{z41}).
\end{PROOF}

\mn
What we really need is the following: it will be used in \ref{c17}(2).

\begin{claim}\label{c14}
Like \emph{\ref{c11}}, but
\begin{enumerate}[$\boxplus_1$]
    \item 
    \begin{enumerate}
        \item [$(a)$] $J \in K_{\clT,2}^\org$ 
\sn        
        \item [$(b)$]  As there. 
\sn
        \item [$(c)$] $X \subseteq J$ is maximal such that 
        $x \neq y \in X$ implies $y \notin \cl_J(\{x\})$\\
        (equivalently, $\neg[x\ E^J\, y]$).\footnote{
            I.e.~we do not demand that $X$ is well-ordered.
        }
\sn
        \item [$(d)$]
        \begin{enumerate}
            \item $Y \in [X]^{\geq\lambda}$ (So in clause $(c)$ we necessarily have $|X| \geq \lambda$.)
\sn
            \item $I \defeq (J \rest \tau_\org) \rest X$ 
\sn
            \item $I \rest \{<\}$ has the $\lambda$-indiscernibility property. (See \emph{\ref{z38}}.)
\sn            
            \item $I \in K_\org$ is $\lambda$-entangled. 
\sn
        \end{enumerate}
\sn
        \item [$(e)$] $Z \defeq$  
        $$
        \big\{F_\bfo^J(a) : a \in X \cap \dom(F_\bfo^J),\ \bfo \in \Omega_a^J,\text{ and } a \in Y \Rightarrow \lh(\bfo) = 1 \mod 2 \big\}
        $$
        \item [$(f)$]
        \begin{enumerate}
            \item If $D \in \clB_\bfp$ then $|Y_D^0| \geq \lambda$, where $Y_D^0 \defeq \{s \in Y : \Omega_s^J = \Omega_D\}$ (recalling Definition \emph{\ref{b37}}).
\sn
            \item $Y$ is the disjoint union of $\LL Y_D^0 : D \in \clB_\bfp\RR$.
        \end{enumerate}
\sn
        \item [$(h)$,] $(i)$,\ \ $(j)$ As there. 
    \end{enumerate}
\sn
    \item $I \rest Y$ is strictly $\mu$-$\Gamma_{\!\oor}$-$\Sigma$-{unembeddable} into $I \rest Z$.
\end{enumerate}
\end{claim}

\sn
\begin{PROOF}{\ref{c14}}
First,
\begin{enumerate}
    \item [$(*)_0$] Let $F$ and $\LL \sigma_s(\bar t_s) : s \in Y\RR$ (where $\bar t_s \in {}^{\omega>}\!Z$) be as in \ref{c8}(A), with $I \rest Y$ and $I \rest Z$ here standing in for $I,J$ there.
\end{enumerate}
It will suffice to find $\bar s_1, \bar s_2 \in {}^2Y$ as in \ref{c8}(B), recalling our choice of $\Gamma$ in clause $\boxplus_1$(b). 

Assume, for the sake of contradiction, that there are no such $\bar s_1, \bar s_2$.
\begin{enumerate}
    \item [$(*)_1$] For {$s \in Y$}, let $(\bar t_s',\bar\bfo_s)$ be such that
    \begin{enumerate}
        \item  $\bar t_s' \in {}^{\omega>}\!X$ (Not to be confused with $\bar t_s \in {}^{\omega>}\!Z$!)
\sn
        \item $\lh(\bar t_s) = \lh(\bar t_s') = \lh(\bar\bfo_s) = n_s = n[s] \defeq \lh(\bar t_s)$ and $\bfo_{s,\ell} \in \Omega_s^J$.
\sn
        \item $\ell < {n_s} \Rightarrow t_{s,\ell} = F_{\bfo_{s,\ell}}^J(t_{s,\ell}')$
\sn
        \item Let 
        $e_s \defeq \big\{ (\ell,k) : \ell,k < {n_s},\ t_{s,\ell}' = t_{s,k}'\big\}.$
\sn
        \item $\LL \bfo_{s,\ell} : \ell < n_s\RR \subseteq \Omega_\bfp$ satisfies 
        $$
        \bfo <_\Omega \bfo_{s,\ell} \Rightarrow s \notin \dom(F_\bfo^J).
        $$
    \end{enumerate}
\end{enumerate}
[Why? For clause (e), recall the definition of $K_{\bfp,2}^\org$ and see \ref{b37}(3)(B). (This is guaranteed when $\bfp$ is well-founded, recalling \ref{b41}.)]
\begin{enumerate}
    \item [$(*)_2$] Choose $\chi$ large enough, and choose $N \prec \big( \clH(\chi),{\in}\big)$ of cardinality $< \lambda$ such that:
    \begin{enumerate}
        \item $J,I,F,\mu,\bfp,X,Y,Z,\Phi,$ and $\LL \sigma_s(\bar t_s) : s \in Y\RR$ all belong to $N$.
\sn
        \item $\|N\| < \lambda$
\sn
        \item $N \cap \lambda \in \lambda$ (so  $\clB,\tau(\Phi)$, $\Omega$, and $\Phi$ are all $\subseteq N$).
    \end{enumerate}
\end{enumerate} 
Next,
\begin{enumerate}
    \item [$(*)_3$] If $s \in Y \setminus N$ then there are sets $v_{s,1}$,  $v_{s,2}$, $v_{s,3}$, and $\clU_s$ (the first three being finite) such that:
    \begin{enumerate}
        \item $v_{s,1} \subseteq (N \cap Z) \cup \{\infty\}$ and 
        $v_{s,2} \subseteq n_s$. 
\sn
        \item $s \in \clU_s \in [Y]^\lambda$ and $\clU_s \cap N = \varnothing$.
\sn
        \item The sequence $\big\LL (\sigma_r,\lh(\bar t_r),e_r,\bar\bfo_r) : r \in \clU_s \big\RR$ is constant.
\sn
        \item  $(\forall\ell \in v_{s,2})(\forall r \in \clU_s)[t_{r,\ell}' \notin N \wedge t_{r,\ell} \notin N ]$ and 
        $$
        \big( \forall\ell \in n_s \setminus v_{s,2} \big) \big( \forall r \in \clU_s \big) \big[t_{r,\ell} = t_{s,\ell} \in v_{s,1} \wedge t_{r,\ell}' = t_{s,\ell}' \big].
        $$
        \item $\LL t_{r,\ell} : \ell < n_s,\ r \in \clU_s,\ t_{s,\ell} \notin N\RR$ and $\LL t_{r,\ell}' : \ell < n_s,\ r \in \clU_s,\ t_{s,\ell}' \notin N\RR$ are without repetition (except for $t_{r,\ell}' = t_{r,k}'$ when $(\ell,k) \in e_s$).
\sn
        \item $t_{s,\ell}' \in v_{s,1}$ for every $\ell \in n_s \setminus v_{s,2}$.
\sn
        \item $v_{s,3} \defeq \{\ell \in v_{s,2} : t_{s,\ell}' = s,\ \ell < n_s\} \subseteq v_{s,2}$ 
    \end{enumerate}
\end{enumerate}        
[Why? Easy.\footnote{
    Earlier, we had also demanded that
    \begin{enumerate}
        \item [(h)] For every $\ell \in v_{s,2}$ there exists $t_\ell^* \in v_{s,1}$ such that if $r \in \clU_s$ \underline{then} $t_\ell^*$ is the $<^J$-minimal member of $\{s \in X \cap N : s \geq_J t_{r,\ell}'\}$. (If there are \emph{no} members of $X \cap N$ above $t_{r,\ell}'$, we say $t_\ell^* \defeq \infty$. This is why we allowed $\infty \in v_{s,1}$ in clause $(*)_1$(a).)
    \end{enumerate}
}]

\begin{enumerate}
    \item [$(*)_4$] Recall that we defined
    $$
    Y_D^0 \defeq \{ s \in Y : \Omega_s^J = \Omega_D\}
    $$
    (so by assumption $\boxplus_1(f)$ we know $Y_D^0 \in [Y]^{\geq\lambda}$).
\sn
    \item [$(*)_5$] For each $D \in \clB$, we can choose $\sigma_D$,  $e_D$, $v_{D,\iota}$ (for $\iota = 2,3$),\\ 
    $\LL t_{D,\ell}^* : \ell \in n_s \setminus v_{s,2}\RR$,
    and $\bar\bfo_D = \LL\bfo_{D,\ell} : \ell < n_s \RR$ such that 
    $|Y_D^1| \geq \lambda$, where
    \begin{align*}
        Y_D^1 \defeq \big\{ s \in Y_D^0 : &\ \sigma_D = \sigma_s,\ \bar\bfo_D = \bar\bfo_s,\ v_{D,\iota} = v_{s,\iota}\text{ for }\iota =2, 3,\\
        &\text{ and $t_{D,\ell}^* = t_{s,\ell}'$ for } \ell \in n_s \setminus v_{s,2} \big\}.
    \end{align*}
\end{enumerate}
[Why does this exist? As $|Y_D^0| = \lambda$ is regular and $> \|N\| \geq |\tau_J| + |\clT_\bfp|$, there exists $s \in Y_D^0 \cap Y \setminus N$, and we can use the same argument as for $(*)_3$.]
\begin{enumerate}
    \item [$(*)_6$]  If $D_1,D_2 \in \clB$, then there exist $\ell = \ell_{D_1,D_2} \in v_{D_1,3}$ and $k = k_{D_1,D_2} \in v_{D_2,3}$ such that $\bfo_{D_1,\ell}$ and $\bfo_{D_2,k}$ have a common $\leq_\Omega$-upper bound. 


\end{enumerate}
Why?  First: 
\begin{enumerate}
    \item [$(*)_{6.1}$] We can choose $s_{1,\eps} \in Y_{D_1}^1$ and $s_{2,\eps} \in Y_{D_2}^1$ for $\eps < \lambda$ such that 
$$
    \eps < \zeta < \lambda \Rightarrow s_{1,\eps} \neq s_{1,\zeta} \wedge s_{2,\eps} \neq s_{2,\zeta}.
$$ 
\end{enumerate}
Second, by $\boxplus_1(d)\,\bullet_3$, without loss of generality 
\begin{enumerate}
    \item [$(*)_{6.2}$] 
    \begin{enumerate}
        \item $\big\LL \LL s_{1,\eps},s_{2,\eps}\RR \caret\, \bar t_{s_{1,\eps}}' \!\caret\ \bar t_{s_{2,\eps}}' : \eps < \lambda \big\RR$ is a qf-indiscernible sequence for the linear order $<^J$. 
\sn
        \item All the individual sequences in clause (a) realize the same quantifier-free type in $I$.
        
        (That is, the type 
        $\tp_\qf(\LL s_{1,\eps},s_{2,\eps}\RR \caret\, \bar t_{s_{1,\eps}}' \!\caret\ \bar t_{s_{2,\eps}}',\varnothing,I)$
        does not depend on $\eps$.)
    \end{enumerate}
    
\end{enumerate}
Third,
\begin{enumerate}
    \item [$(*)_{6.3}$] $\big\LL \LL s_{1,\eps},s_{2,\eps}\RR \caret\,  \bar t_{s_{1,\eps}} \!\caret\ \bar t_{s_{2,\eps}} : \eps < \lambda \big\RR$ is a qf-indiscernible sequence for $<^J$ as well.
\end{enumerate}
[Why? By $(*)_{6.2}$ and $\boxplus_1(i),(j)$, recalling $(*)_5$.]

\mn
That is,
\begin{enumerate}
    \item [$(*)_{6.4}$]
    \begin{enumerate}[$\bullet_1$]
        \item For $\eps < \zeta < \lambda$, let $\clW_{\eps,\zeta}^0 \defeq \big\{(\ell,k) \in n_{D_1} \times n_{D_2} : t_{s_{1,\eps},\ell}'\ R^J\, t_{s_{2,\zeta},k}' \big\}$.
\sn
        \item $\clW_* \defeq \big\{(\ell,k) \in n_{D_1} \times n_{D_2} : \bfo_{D_1,\ell} = \bfo_{D_2,k} \big\}$
    \end{enumerate}
\sn
    \item [$(*)_{6.5}$] If $\eps_1 < \zeta_1 < \lambda$ and $\eps_2 < \zeta_2 < \lambda$ are such that $\clW_{\eps_1,\zeta_1}^0 = \clW_{\eps_2,\zeta_2}^0$, \underline{then}
    \begin{enumerate}[$\bullet_1$]
        \item The sequences $\LL s_{1,\eps_1},s_{2,\eps_1}\RR \caret \bar t_{s_{1,\eps_1}} \!\caret\, \bar t_{s_{2,\zeta_1}}$ and 
        $\LL s_{1,\eps_2},s_{2,\eps_2}\RR \caret \bar t_{s_{1,\eps_2}} \!\caret\, \bar t_{s_{2,\zeta_2}}$ realize the same quantifier-free type in $J \rest \{<\}$.
\sn
        \item The sequences $\bar t_{s_{1,\eps_1}} \!\caret\, \bar t_{s_{2,\zeta_1}}$ and 
        $\bar t_{s_{1,\eps_2}} \!\caret\, \bar t_{s_{2,\zeta_2}}$ realize the same quantifier-free type in $J$.
\sn
        \item If 
        $$
        [s_{1,\eps_1}\, R^J\, s_{1,\zeta_1}] \Leftrightarrow \neg[ s_{2,\eps_1}\, R^J\, s_{2,\zeta_1}]
        $$
        \underline{then} we get the desired {contradiction}.
    \end{enumerate}
\end{enumerate}
[Why? For $\bullet_1$, first recall $(*)_{6.2}$(a) for the order on 
$\LL s_{1,\eps},s_{2,\eps}\RR \caret \LL \bar t_{s_{1,\eps}}' \!\caret\ \bar t_{s_{2,\eps}}' : \eps < \lambda\RR$. For $\bullet_2$, recall $(*)_{6.2}$(b) 
and our assumption $\clW_{\eps_1,\zeta_1}^0 = \clW_{\eps_2,\zeta_2}^0$
for the graph relation. Lastly, $\bullet_3$ follows from our choices.]

\medskip
\begin{enumerate}
    \item [$(*)_{6.6}$] $(\forall \eps < \zeta < \lambda) \big[ \clW_{\eps,\zeta}^0 \subseteq \clW_* \big]$ 
\end{enumerate}
[Why? Holds by assumption $\boxplus_1(h)$.]

\medskip
\begin{enumerate}
    \item [$(*)_{6.7}$] If $\clW_*$ and $v_{D_1,3} \times v_{D_2,3}$ are disjoint, then we get a contradiction.
\end{enumerate}
[Why? Because $I \in K_\org$ is $\lambda$-entangled,
and get a contradiction by $(*)_{6.5}$+$(*)_{6.6}$.]

\medskip
So for some $\bfo \in \Omega$ and $\ell < n_{D_1}$, $k < n_{D_2}$, for {all} $\eps < \zeta < \lambda$, we have $F_\bfo(s_{1,\eps}) = t_{s_{1,\eps},\ell}$ and $F_\bfo(s_{2,\zeta}) = t_{s_{2,\zeta},k}$. 

Now $J \in K_{\clT,2}^\org \subseteq K_{\clT,1}^\org$ (see \ref{b19}(2),(3)). Therefore, recalling $(*)_1$(e), we conclude that $\bfo_{D_1,\ell}$ and $\bfo_{D_2,k}$ have a common $\leq_\Omega$-upper bound. 

This proves $(*)_6$.

\mn
Now we recall \ref{b3}(1A):
\begin{enumerate}
    \item [$\boxplus$] $\theta_\bfp > \aleph_0$, and the intersection of countably many members of $\clB$ will always contain some other member of $\clB$.
\end{enumerate}
Hence
\begin{enumerate}
    \item [$(*)_7$] There exists $\bfn$ such that 
    $$
    \clB_\bfn \defeq \{D \in \clB : n_D = \bfn\}
    $$ 
    is $\leq_\bfp$-cofinal in $(\clB,\leq_\bfp)$.
\sn
    \item [$(*)_8$] Let $\overbar\bfm_D \defeq \big(\LL \lh(\bfo_{D,\ell}) : \ell < \bfn\RR, v_{D,2}, v_{D,3} \big)$.
\sn
    \item [$(*)_9$] Let $\bbE$ be an ultrafilter on $\clB$ which includes $\clB_\bfn$, such that for every $D \in \clB$ we have 
    $$
    \{D' \in \clB : D' \subseteq D\} \in \bbE.
    $$
\end{enumerate}
[Why does such an $\bbE$ exist? Because $(\clB_\bfp,\supseteq)$ is directed and $\clB_\bfn$ is cofinal in it.]

\sn
\begin{enumerate}
    \item [$(*)_{10}$] 
    \begin{enumerate}
        \item For each $D \in \clB_\bfn$ there exist $\ell_D,k_D < \bfn$ such that
        $$
        \clX_D \defeq \{D' \in \clB_\bfn : \ell_D = \ell_{D,D'},\ k_D = k_{D,D'}\} \in \bbE.
        $$
        \item There exist $\ell_*,k_* < \bfn$ such that
        $$
        \clB_\bullet \defeq \{D \in \clB_\bfn : \ell_D = \ell_*,\ k_D = k_*\} \in \bbE.
        $$
    \end{enumerate}
\end{enumerate}
[Why? Obvious: there are only finitely many possibilities (clause (a) gives $<\bfn$ and clause (b) gives $<\bfn^2$), and $\bbE$ is an ultrafilter. We may add more, but this is not necessary.]

\mn
\begin{enumerate}
    \item [$(*)_{11}$] $\big(\forall^\bbE D_1 \in \clB_\bullet\big) \big(\forall^\bbE D_2 \in \clB_\bullet\big) \big[ \lh(\bfo_{D_1,\ell_*}) \leq \lh(\bfo_{D_2,k_*}) \big]$
\end{enumerate}
[Why? By clause (C)(b) of Definition \ref{b3}(1), for any $D_1 \in \clB_\bullet$ there is $D_2' \in \clB_\bfp$ such that
$$
\{\eta :   
\eta  
\text{ appears in } \bfo_{D_1,\ell_*} \} \cap D_2' = \varnothing.
$$
(That is, $\eta \in  
\rang(\bar\eta_{D_1,\ell_*})$  
where $\bfo_{D_1,\ell_*} = (\bar\eta,\bar\iota) = (\bar\eta_{D_1,\ell_*},\bar\iota_{D_1,\ell_*})$.)

Now any $D_2 \subseteq D_2'$ from $\clB_\bullet$ will work.]

\mn
\begin{enumerate}
    \item [$(*)_{12}$] If $k < \omega$ and  $D_i \in \clB_\bullet$ for $i < k$, \underline{then} $\{\bfo_{D_i,\ell_*} : i < k\}$ has a common $\leq_\Omega$-upper bound.
\end{enumerate}
[Why? Recalling $(*)_{10}$(a)+$(*)_{11}$, choose {any} $D_* \in \bigcap\limits_{i <k}\clX_{D_i}$ such that $\lh(\bfo_{D_*,k_*}) \geq \lh(\bfo_{D_i,\ell_*})$ for all $i < k$. Then $\bfo_{D_*,k_*}$ is a common upper bound.]

\mn
\begin{enumerate}
    \item [$(*)_{13}$] Assume $D \in \clB_\bullet$, and let $\eta_D \in \clT$ denote the 
    $0^\tthh$ $\eta$-term in $\bfo_{D,\ell_*}$.
    \begin{enumerate}
        \item $\bfG \defeq \{\eta_D : D \in \clB_\bullet\}$ is directed under $\leq_\clT$.
\sn
        \item $\bfG \cap D \neq \varnothing$ for all $D \in \clB_\bullet$.
    \end{enumerate}
\end{enumerate}
[Why? Clause (a) holds by our choices. Clause (b) holds because for every $D \in \clB_\bfp$ there exists $D' \in \clB_\bullet$ such that $D' \subseteq D$.]

\mn
So $\bfG$ contradicts clause \ref{b3}(1A)(d) in the definition of `strong twinship parameter,' and we are done.
\end{PROOF}

\mn
Now at last we are able to prove Theorem \ref{x2}.
\begin{conclusion}\label{c17}
Assume $\bbP$ is a forcing notion adding no new $\omega$-sequence of ordinals, but does add some {infinite} sequence (so necessarily of length $\geq \omega_1$).

Assume $T \subseteq T_1$ are complete first-order theories, and $T$ is independent as witnessed by $\varphi = \varphi(\bar x_{[k]},\bar y_{[k]}) \in \bbL(\tau_T)$. 

\mn
$1)$ \underline{Then} there are models $M,N$ such that:
\begin{enumerate}[$(a)$]
    \item $M$ and $N$ are models of $T_1$ of cardinality $\lambda \defeq (2^{\|\bbP\|} + |T_1|)^+$ (or use any regular cardinal $\lambda' \ge \lambda$).
\sn
    \item $\Vdash_\bbP ``M\rest \tau_T \cong N\rest \tau_T"$, and even $\Vdash_\bbP ``M \cong N"$.
\sn
    \item $M\rest \tau_T$ and $N \rest \tau_T$ are not isomorphic.
\sn
    \item Moreover, $M$ is $(\lambda,2^{\|\bbP\|},\varphi)$-\emph{far} from $N$ (see \emph{\ref{a68}(1)}). 
\end{enumerate}

\mn
$2)$ We may strengthen clause $(d)$ to 
\begin{enumerate}
    \item[$(d)^+$] $M$ and $N$ are $(\lambda,2^{\|\bbP\|},\varphi)$-\emph{far} from each other.
\end{enumerate}
\end{conclusion}

\mn
\begin{discussion}\label{c3}
Recalling Definition \ref{a59}, note that in \ref{c17} we cannot deduce that $M$ and $N$ are $(\lambda,\varphi)$-far, as the partial isomorphisms $F_{\eta,\iota}^J$ form a witness. Now clause $\boxplus_1$(c) in the assumptions of \ref{c11} is strong, \emph{but} it is only talking about $(J \rest \tau_\org) \rest X$, so the partial isomorphism $F_{\eta,\iota}^J$ disappears.

However, the possibility of being $(\lambda,|\clB_\bfp|,\varphi)$-far (see Definition \ref{a68}) is not excluded.
\end{discussion}

\sn
\begin{remark}\label{c21}
The following are natural extensions of Claim \ref{c17}. (Their proofs will be delayed to \cite{Sh:F2433}.)

\mn
$3)$ In parts $(1)$ and $(2)$ of \ref{c17}, we may add 
\begin{enumerate}
    \item [$(e)$] $M$ and $N$ are $\bbL_{\infty,\lambda}$-equivalent.
\end{enumerate}

\mn
$4)$ If $\lambda = \cf(\lambda) > 2^{\|\bbP\|} + |T_1|$ and $\xi < \lambda$, \underline{then} we can find models $M,N$ of $T_1$ such that:
\begin{enumerate}[$(a)$\,]
    \item $\|M\| = \|N\| = |{}^{\xi>}\!\lambda|$
\sn
    \item[$(b),$] $(c)$ \  As in \ref{c17}(1).
\sn
    \item [$(d)$\,] $M$ and $N$ are cofinally $(\lambda,\xi)$-equivalent (see Definition {\ref{z32}}).
\sn
    \item [$(e)$\,] $M$ and $N$ are $(\lambda,2^{\|\bbP\|},\varphi)$-\emph{far} from each other.
\end{enumerate}

\mn
5) It is enough to demand that $\lambda$ is regular, $>|\bbP|$, and $\geq$ the number of maximal antichains in $\bbP$.
\end{remark}

\sn
\begin{PROOF}{\ref{c17}}
\textsc{Proof of \ref{c17}:}

\sn
1) First, choose a strong twinship parameter $\bfp$ by Claim \ref{b62} (so $\theta_\bfp \leq \|\bbP\|$,\\ $|\clT_\bfp| \leq \|\bbP\|$, and $|\clB_\bfp| \leq 2^{\|\bbP\|}$).
\begin{enumerate}
    \item [$(*)_1$]  We may require 
    \begin{enumerate}
        \item $\theta_\bfp \defeq \min\{\theta :\ \Vdash_\bbP\text{``there is a new } \eta \in {}^\theta2"\}$
\sn
        \item $|\clB_\bfp| = \big|\{ |Y| : Y \text{ is a maximal antichain of } \bbP\} \big|$.
    \end{enumerate}
\end{enumerate}
Second,
\begin{enumerate}
    \item [$(*)_2$] 
    \begin{enumerate}
        \item Choose $\lambda$ regular such that $\lambda > |\tau(T_1)| + |\clT_\bfp|$, $\lambda \geq |\clB_\bfp|$.
\sn
        \item $\bfc : [\lambda]^2 \to \omega$ witnesses
        the property called 
        $\Pr_0(\lambda,\aleph_0)$ in \cite[Th.\,1.1]{Sh:280} (later called $\Pr_0(\lambda,\lambda,\aleph_0,\aleph_0)$).
    \end{enumerate}
\end{enumerate}
Recall
\begin{enumerate}
    \item[$\odot$] $\Pr_0(\lambda,\aleph_0)$ means that if $n < \omega$ and $\big\LL \bar\zeta_\alpha = \LL\zeta_1^\alpha : i \leq n \RR : \alpha < \lambda \big\RR$ is such that $\zeta_1^\alpha < \zeta_2^\alpha < \ldots < \zeta_n^\alpha < \lambda$ and $\alpha < \beta < \lambda \Rightarrow \bar\zeta_\alpha \cap \bar\zeta_\beta = \varnothing$, \underline{then} for any function $h : n \times n \to \mu$ there exists $\alpha < \beta < \lambda$ such that $\zeta_n^\alpha < \zeta_1^\beta$ and 
    $$
    k,\ell \in [1,n] \Rightarrow \bfc(\{\zeta_k^\alpha,\zeta_\ell^\beta\}) = h(k,\ell).
    $$
\end{enumerate}

\sn
[Why {does such a $\bfc$ exist}? By \cite[Th. 1.1]{Sh:280}.]

\mn
Now choose
$J \in K_{\clT,2}^\org$ as in \ref{c11}\,$\boxplus_1$ such that:
\begin{enumerate}
    \item [$(*)_3$]
    \begin{enumerate}
        \item $X \subseteq J$, $(X,<_J) = (\lambda,<)$, and the pair $(X,J)$ satisfies clauses (h),(i),(j) of \ref{c11}\,$\boxplus_1$.
\sn
        \item For $s,t \in X$, we have $(s,t) \in R^J \Leftrightarrow s,t \in J \wedge s \neq t \wedge \bfc(\{s,t\}) = 1$.
    \end{enumerate}
\end{enumerate}
[Why? For each $D \in \clB_\bfp$ we can find $I_D \in K_\org$ of cardinality $|D|$ (which is infinite, but $\leq \aleph_0 + |\clT_\bfp| < \lambda$) such that {$s \in I_D \Rightarrow  \Omega_s^{I_D} = \Omega_D$} by {\ref{b34}}. Let $s_D = s[D]$ be some member of $I_D$.

Let $\LL D_\alpha : \alpha < \lambda\RR \in {}^\lambda\clB_\bfp$ be such that $(\forall D \in \clB_\bfp)(\exists^\lambda \alpha < \lambda)[D_\alpha = D]$. We define $J \in K_{\clT,2}^\org$ as follows.
\begin{enumerate}
    \item[$(*)_4$] 
    \begin{enumerate}
        \item $|J| \defeq \lambda \times I_D$. (That is, $\{(\alpha,s) : \alpha < \lambda,\ s \in I_D\}$.)
\sn
        \item $(\alpha_1,s_1) <_J (\alpha_2,s_2)$ \underline{iff} 
        \begin{enumerate}
            \item $(\alpha_1,s_1), (\alpha_2,s_2) \in J$
\sn
            \item $\alpha_1 < \alpha_2 \vee \big[\alpha_1 = \alpha_2 \wedge s_1 <_{I_{D_\alpha}} s_2 \big]$.
        \end{enumerate}
\sn
        \item $F_{\eta,1}^J(\alpha_1,s_1) = (\alpha_2,s_2)$ \underline{iff} $\alpha_1 = \alpha_2 \wedge D = D_{\alpha_1} \wedge F_{\eta,1}^{I_D}(s_1) = s_2$.
\sn
        \item $X \defeq \big\{ (\alpha,s_{D_\alpha}) : \alpha < \lambda\big\}$
\sn
        \item We choose 
        \begin{align*}
            R^J \defeq \big\{ (F_\bfo^J(\alpha_1,s_1),F_\bfo^J(\alpha_2,s_2)) : &\ \alpha_1 \neq \alpha_2 < \lambda,\\
            &\ \bfc(\{\alpha_1,\alpha_2\}) = 1,\\
            &\ \bfo \in \Omega_{D_{\alpha_1}} \cap \Omega_{D_{\alpha_2}}\big\}.
        \end{align*}
    \end{enumerate}
\end{enumerate}
Now check. Note that $(J,<_J)$ is not $(\lambda,<)$, {but we} get this by renaming.]

\medskip
Next, 
\begin{enumerate}
    \item [$(*)_5$]
    \begin{enumerate}
        \item Choose $\varphi = \varphi(\bar x_{[k]},\bar y_{[k]}) \in \bbL_{\tau_T}$ witnessing the independence property for $T$ (see Definition \ref{a44}(2)).
\sn
        \item Choose $\Phi \in \Upsilon_{K_\org}\big[T_1,|T_1|\big]$ such that
$$
        I \in K_\org \Rightarrow \big( \forall s,t \in I\big) \big[ \GEM(I,\Phi) \models ``\varphi[\bar a_s,\bar a_t]\ \Leftrightarrow\ s\ R^I\, t"\big].
$$
    \end{enumerate}
\end{enumerate}
[Why can we do this? By \ref{a51}.]

\medskip
Now we shall finish by \ref{c11}. That is, for $\ell = 1,2$:
\begin{enumerate}
    \item [$(*)_6$]
    \begin{enumerate}
        \item Recall that we defined 
$$
        X_\ell \defeq \{F_\bfo^J(s) : s \in X \cap \dom(F_\bfo^J),\, \bfo \in \Omega_\bfp^\fr,\, \lh(\bfo) \equiv \ell \mod 2\}.
$$
        \item Now let $M_\ell$ be the submodel of $\GEM(J \rest \tau_\org,\Phi)$ generated by\\ $\{a_s : s \in X_\ell\}$, where $\LL a_s : s \in J\RR$ is the skeleton.
    \end{enumerate}
\end{enumerate}
Now $(M,N) \defeq \big( M_1, M_2\big)$ are as required. Moreover, $M_1$ and $M_2$ are strictly $\bfp$-isomorphic (see \ref{b9}).

We have to check all four clauses of the conclusion in part (1). 

\mn
\textbf{Clause $(a)$:} $M$ and $N$ are models of $T_1$ of cardinality $\lambda$.

Why? ``Models of $T_1$:" this follows from our choice of $\Phi$ in $(*)_5(b)$ and of $M_\ell$ in $(*)_6(b)$. Their cardinalities are as required because they are the same as $|X_\ell|$, recalling $(*)_6(a)$.

\bn
\textbf{Clause $(b)$:} $\Vdash_\bbP ``M\rest \tau_T \cong N\rest \tau_T"$, and even $\Vdash_\bbP ``M \cong N"$.

Why? This follows from 
$$\Vdash_\bbP ``(J \rest \tau_\org) \rest X_1 \cong (J \rest \tau_\org) \rest X_2",
$$
which holds by \ref{b47}(1).

\bn
\textbf{Clause $(c)$:} $M\rest \tau_T$ and $N \rest \tau_T$ are not isomorphic.

This is the main point, and it follows from \ref{c1}.

\bn
\textbf{Clause $(d)$:} $M$ is $(\lambda,2^{\|\bbP\|},\varphi)$-{far} from $N$.

Just repeat the proof of \ref{c14}.

\bn
2) We need to prove
\begin{enumerate}
    \item [$(D)^+$] $M$ and $N$ are $(\lambda,2^{\|\bbP\|},\varphi)$-\emph{far} from each other.
\end{enumerate}

We intend to use \ref{c14} instead of \ref{c11} for {$\ell = 1,2$}.
For proving clause $(d)^+$, we need $Y_1$ and $Y_2$, hence $Z_1,Z_2$.

So let $Y_1,Y_2 \in [X]^\lambda$ be a partition of $X$ and let 
$$
Z_\ell \defeq \{F_\bfo^J(s) : s \in X \cap \dom(F_\bfo^J), \text{ and } s \in Y_\ell \Rightarrow \lh(\bfo) \equiv \ell \mod 2\}.
$$
Also choose $\eta_s \in D_s^J$ for all $s \in X$. Lastly, let 
$$
X_\ell \defeq \{F_\bfo^J(s) : s \in X \cap \dom(F_\bfo^J),\, \lh(\bfo) \equiv \ell \mod 2\}.
$$
We continue as in the proof of part (1).
\end{PROOF}

\mn
\begin{claim}\label{c23}
Like \emph{\ref{c17}}, when $T$ is countable and superstable, with OTOP or DOP, and $T_1 = T$.
\end{claim}

\sn
\begin{remark}\label{c29}
1) For the existence of ccc forcings with $T$ unsuperstable, this is proved in  \cite[1.1]{Sh:464}, building on \cite[Ch.X, \S2; Ch.XIII, \S2]{Sh:c}. Here we have to use $\Phi \in \Upsilon[T,K_\org]$ (not just $\Upsilon[T,K_\oor]$).

\mn
2) For $T$ stable but not superstable, there is a more liberal version allowing infinite sequences in the parameters; we will return to this in \cite{Sh:F2433}.

\mn
3) Concerning OTOP below, note that \cite[1.28\subref{d17}]{Sh:E59} is quite relevant.

\mn
4) This proof will require some knowledge of stability theory.
\end{remark}

\sn
\begin{PROOF}{\ref{c23}}
Similar to the proof of \ref{c17}.

The point is that those properties imply the existence of $\Phi$ as in the proof of \ref{c17}, except that the relevant $\varphi$ is not first-order.  

Below, we will prove existence for the two {possible} cases.

\bn
\textbf{Case 1: $T$ has $\mathrm{OTOP}$ but fails $\mathrm{DOP}$.}

Recall the definition of OTOP from \cite[Ch.XII, 4.1, p.608]{Sh:c}.
\begin{enumerate}
    \item [$\boxplus_{1.1}$] We say that $T$ has the \emph{Omitting Type Order Property} (OTOP) \underline{if} there is a type $p(\bar x,\bar y,\bar z)$ (where $\bar x,\bar y,\bar z$ are finite sequences with $\lh(\bar x) = \lh(\bar y)$) such that for every $\lambda$ and every two-place relation $R$ on $\lambda$, there exists a model $M$ of $T$ and $\LL\bar a_\alpha : \alpha < \lambda\RR \subseteq M$ such that
    $$
    \alpha < \beta < \lambda \Rightarrow \big[ \alpha\ R\ \beta \Leftrightarrow \text{The type $p(\bar a_\alpha,\bar a_\beta,\bar z)$ is realized in } M\big].
    $$
\end{enumerate}

\mn
Everything we need has been proven in \cite{Sh:E59}, but let us elaborate. 

For Case 1, the following claim will suffice.
(Note that this is \cite[1.28\subref{d17}]{Sh:E59}, for ordered graphs instead of linear orders.)

\begin{enumerate}
    \item [$\boxplus_{1.2}$] If $T$ is first-order, countable, and has the OTOP, \underline{then} for some sequence $\olsi\varphi = \LL \varphi_i(\bar x,\bar y,\bar z) : i < i_*\RR$ of first-order formulas in $\bbL(\tau_T)$ and a template $\Phi$ proper for ordered graphs, we have:
    \begin{enumerate}
        \item $\tau_T \subseteq \tau_\Phi$ and $|\tau_\Phi| = |\tau_T|+ \aleph_0$.
\sn
        \item $\GEM_{\tau_T}(I,\Phi) \models T$ for all $I \in K_\org$.
\sn
        \item If $I \in K_\org$ and $s,t \in I$, \underline{then}
        $$
        \big[ \GEM_{\tau_T}(I,\Phi) \models (\exists \bar x) \textstyle\bigwedge\limits_{i<i_*} \varphi_i(\bar x,\bar a_s,\bar a_t) \big] \Leftrightarrow I \models `s\ R\ t\text{'}.
        $$
    \end{enumerate} 
\end{enumerate}

\begin{PROOF}{\boxplus_{1.2}}
\textsc{Proof of $\boxplus_{1.2}$:}

We would like to apply \cite[1.25(e)\subref{d8}]{Sh:E59}, but it requires us to assume enough cases of a certain partition theorem generalizing Erd\H{o}s-Rado (with $K_\oor$ replaced by $K_\org$). However, this theorem is proven to hold only after a forcing --- in fact, for every {infinite cardinal} $\kappa$ there is a $\kappa$-complete class forcing which ensures that the analogous result holds above $\kappa$ (it holds by \cite{Sh:289}).

Moreover, by \cite[Ch.XII, \S5]{Sh:c}, assuming $T$ is countable, superstable, and has the NDOP, \underline{if} it has the $(\aleph_0,2)$-existence property \underline{then} it satisfies clause (B) or (C) of \ref{x17}(1)$\boxplus_1$. {This gives us a} contradiction, so $T$ fails the $(\aleph_0,2)$-existence property. This is preserved by any 
$\kappa$-complete forcing $\bbP$, even for $\kappa \defeq \aleph_0$ (but $\aleph_1$ is more convenient here).

By \cite[Ch.XII, 4.3, p.609]{Sh:c} $T$ still has the OTOP in 
$\bfV^\bbP$, so we can apply \cite[1.25(e)\subref{d8}]{Sh:E59} to get 
$\Phi$ as promised. But as {we assumed} $\kappa > \aleph_0$ we know 
$\Phi \in \bfV$, so we are done proving $\boxplus_{1.2}$, and hence we have finished {the present case}.
\end{PROOF}

\bn
\textbf{Case 2: $T$ has} DOP \textbf{(}\emph{dimensional order property}\textbf{).}

Here we shall use \cite[Ch.X, 2.1-2, p.512]{Sh:c} and \cite[Ch.X, 2.4, p.515]{Sh:c}. Without loss of generality $\kappa \defeq |T|^+$.

\sn
Recall:
\begin{enumerate}
    \item  [$\boxplus_{2.1}$] We say that $T$ has the \emph{Dimensional Order Property} (DOP) \underline{if} there are $\kappa$-saturated models $M_\ell \prec \gC$ of cardinality $\leq 2^{|T|}$ (for $\ell = 0,1,2$), where $\gC$ is a $\varkappa$-saturated model of $T$ for some $\varkappa > \|M_1\| + \|M_2\|$, such that $M_0 \subseteq M_1 \cap M_2$, $\{M_1,M_2\}$ is independent over $M_0$, and the $\kappa$-prime model $M$ over 
    $M_1 \cup M_2$ is not $\bfF_\kappa^a$-minimal over $M_1 \cup M_2$.
\sn
    \item  [$\boxplus_{2.2}$] For any $\kappa > 2^{|T|}$ (for transparency), for any $\kappa$-complete forcing $\bbP$, the relevant properties of $T$ are still preserved in $\bfV^\bbP$.
\sn
    \item  [$\boxplus_{2.3}$] In $\boxplus_{2.1}$, we can find finite $\bar a_0,\bar a_1,\bar a_2, \bar b,c$ such that:
    \begin{enumerate}
        \item $\bar a_\ell \subseteq M_\ell$ for $\ell =0,1,2$.
\sn
        \item $\bar b \caret \LL c\RR \subseteq M$
\sn
        \item For $\ell = 1,2$, $\tp(\bar a_\ell, M_0,M)$ does not fork over {$\bar a_0$} and $\tp(\bar a_\ell,\bar a_0,M)$ is stationary.
\sn
        \item $\tp(\bar b \caret \LL c\RR, M_1 \cup M_2,M)$ does not fork over $\bar a_1 \caret \bar a_2$, and 
        $$
        \tp(\bar b \caret \LL c\RR, \bar a_1 \caret \bar a_2,M)
        $$ 
        has a unique non-forking extension in $\clS(M_1,M_2)$.
\sn
        \item $\tp(c,\bar b,M)$ is stationary {and} $c \in \bfI$, {where} $\bfI \subseteq M$ is infinite, indiscernible over 
        $M_1 \cup M_2 \cup \bar b$, and based on $\tp(c,\bar b,M)$.
    \end{enumerate}
\end{enumerate}
[Why? Because $T$ is superstable. (For clause (d), recall \cite[Ch.XII, \S3]{Sh:c}.)]
\begin{enumerate}
    \item  [$\boxplus_{2.4}$] For every $\lambda$ and $R \subseteq \lambda \times \lambda$, we can find $\LL f_\alpha^\ell : \alpha < \lambda,\ \ell = 1,2\RR$,\\ $\LL f_{\alpha,\beta}, \bfI_{\alpha,\beta} : (\alpha,\beta) \in R\RR$, and $N$ such that:
    \begin{enumerate}
        \item $M_0 \prec N$
\sn
        \item $f_\alpha^\ell$ is an elementary embedding of $M_\ell$ into $N$ over $M$.

        We let $\bar a_{\ell,\alpha} \defeq f_\alpha^\ell(\bar a)$ and $M_{\ell,\alpha} \defeq f_\alpha^\ell(M_\ell)$.
\sn
        \item $\LL M_{\ell,\alpha} : \alpha < \lambda,\ \ell = 1,2\RR$ is independent over $M_0$.
\sn
        \item $f_{\alpha,\beta}$ is an elementary embedding of $M$ into $N$ extending $f_\alpha^1 \cup f_\beta^2$. 

         We let $\bar b_{\alpha,\beta} \defeq f_{\alpha,\beta}(\bar b)$ and $c_{\alpha,\beta} \defeq f_{\alpha,\beta}(c)$.
\sn
        \item $\bfI_{\alpha,\beta}$ is an {indiscernible} set over $\bar b_\alpha$ based on $\tp(c_{\alpha,\beta},\bar b_\alpha)$ (equivalently, on $\tp(c_{\alpha,\beta},\bar b_\alpha \cup M_{1,\alpha} \cup M_{2,\beta})$) of cardinality $|T|$.
\sn
        \item $N$ is $|T|^+$-prime over 
        $$
        \bigcup\limits_{\substack{\alpha < \lambda\\\ell = 1,2}} M_{\ell,\alpha} \cup \bigcup\limits_{(\alpha,\beta) \in R} f_{\alpha,\beta}(\bar b \caret \LL c\RR) \cup \bigcup\limits_{(\alpha,\beta) \in R} \bfI_{\alpha,\beta}.
        $$
        \item If $(\alpha,\beta) \in \lambda\times \lambda \setminus R$, \underline{then} we cannot find $\bar b_{\alpha,\beta}, c_{\alpha,\beta},\bfI_{\alpha,\beta}$ as above (in $N$).
    \end{enumerate}
\end{enumerate}
As in Case 1, we can find a suitable $\Phi$.
\end{PROOF}

\mn
\begin{problem}\label{c26}
Can we eliminate ``$\theta_\bfp > \aleph_0$" in \S3?
\end{problem}

\bibliographystyle{amsalpha}
\bibliography{shlhetal}
\end{document}